\documentclass{amsart}
\title{Sufficient and Necessary Conditions for Eckart-Young-like Result for Tubal Tensors}

 \author{Uria Mor}
 \thanks{School of Mathematical Sciences, Raymond and Beverly Sackler Faculty of Exact Sciences, Tel Aviv University, Tel Aviv 6997801, Israel.}
\date{\today}

\usepackage[dvipsnames]{xcolor}
\usepackage{amssymb,amsmath, amsthm, amsfonts,comment,algorithm}
\usepackage[margin=1in]{geometry}
\usepackage{algpseudocode}
\usepackage{algorithmicx}
\usepackage{mathtools}
\usepackage{amssymb}
\usepackage{mathrsfs}

\usepackage[colorlinks=true]{hyperref}
\usepackage[nameinlink]{cleveref}

\crefname{subsection}{subsection}{subsections}
\crefname{subsubsection}{subsubsection}{subsubsections}
\newtheorem{theorem}{Theorem}[section]
\newtheorem{lemma}[theorem]{Lemma}

\newtheorem{proposition}[theorem]{Proposition}
\newtheorem{definition}[theorem]{Definition}

\crefname{req}{Requirement}{requirements}
\Crefname{req}{Requirement}{Requirements}
\crefname{prprt}{Property}{properties}
\Crefname{prprt}{Property}{Properties}
\crefname{def}{Definition}{definitions}
\Crefname{def}{Definition}{Definitions}
\crefname{claim}{Claim}{claims}
\Crefname{claim}{Claim}{Claims}
\crefname{q}{question}{questions}
\crefname{lemma}{Lemma}{lemmas}
\Crefname{lemma}{Lemma}{Lemmas}
\crefname{obs}{Observation}{observations}
\Crefname{obs}{Observation}{Observations}
\crefname{prop}{Proposition}{proposition}
\Crefname{prop}{Proposition}{Propositions}
\crefname{cor}{Corollary}{corollaries}
\Crefname{cor}{Corollary}{Corollaries}
\crefname{equation}{Eq.}{Eqs.}
\Crefname{equation}{Eq.}{Eqs.}

\crefname{example}{example}{examples}
\Crefname{example}{Example}{Examples}

\crefname{conj}{conjecture}{conjectures}
\Crefname{conj}{Conjecture}{Conjectures}

\providecommand{\examplename}{Example}

\newtheorem{remark}[theorem]{Remark}

\usepackage[h]{esvect}

\RenewDocumentCommand{\Vec}{me{_}}{%
  \IfNoValueTF{#2}{\vv{#1}}{\vv*{#1}{\mspace{-2mu}#2}}%
}

\usepackage{multibib}
\newcites{latex}{\LaTeX-Literature}
\begin{document}

\newcommand{\TRN}{\textnormal{train}}
\newcommand{\TST}{\textnormal{test}}
\newcommand{\tn}[1]{\textnormal{#1}}

\newcommand{\wh}[1]{\widehat{#1}{}}
\newcommand*\wt[1]{\widetilde{#1}{\.}}
\newcommand*\wol[1]{\overline{\vphantom{#1^{a}}{#1}}}

\newcommand{\noun}[1]{\textsc{#1}\xspace}
\newcommand{\bs}[1]{\boldsymbol{#1}}
\newcommand{\hbs}[1]{\ensuremath{\wh{\bs{#1}}}}
\newcommand{\rank}{\text{\ensuremath{\operatorname{rank}}}}
\newcommand{\argmin}{\textnormal{argmin}}
\newcommand{\argmax}{\text{\ensuremath{{\arg}\max}}}

\newcommand{\xx}{\times}
\newcommand{\code}[1]{\codex{#1}}
\newcommand{\pinv}{\ensuremath{\dagger}}
\renewcommand{\P}{\pinv}
\newcommand{\mpn}{\ensuremath{m \xx p \xx n}}
\newcommand{\pmn}{\ensuremath{p \xx m \xx n}}
\newcommand{\innerprod}[2]{\left\langle#1,#2 \right\rangle}
\newcommand{\mpd}{\ensuremath{m \xx p \xx \ddim}}
\newcommand{\pmd}{\ensuremath{p \xx m \xx \ddim}}
\newcommand{\dimsI}{\ensuremath{I_1 \xx \cdots \xx I_d}}
\newcommand{\dimsr}{\ensuremath{r_1 \xx \cdots \xx r_d}}
\newcommand{\scalemath}[2]{\scalebox{#1}{\mbox{\ensuremath{\displaystyle #2}}}}
\newcommand{\dbr}[1]{\ensuremath{\llbracket #1 \rrbracket}}

\newcommand{\tens}[1]{\bs{\mathcal{#1}}}
\newcommand{\tA}{\tens{A}}
\newcommand{\tAt}{\tA^{\T}}

\newcommand{\tT}{\tens{T}}

\newcommand{\btA}{\bar{\tens{A}}}
\newcommand{\ttA}{\tilde{\tens{A}}}
\newcommand{\ttAt}{\ttA^{\T}}
\newcommand{\thA}{\wh{\tA}}
\newcommand{\tthA}{\wh{\ttA}}
\newcommand{\thAt}{\thA^{\T}}
\newcommand{\tthAt}{\tthA^{\T}}
\newcommand{\tB}{\tens{B}}
\newcommand{\tBt}{\tB^{\T}}
\newcommand{\thB}{\wh{\tB}}
\newcommand{\thBt}{\thB^{\T}}
\newcommand{\tR}{\tens{R}}
\newcommand{\thR}{\wh{\tR}}
\newcommand{\tRt}{\tR^{\T}}
\newcommand{\tM}{\tens{M}}
\newcommand{\thM}{\wh{\tM}}
\newcommand{\tMt}{\tM^{\T}}

\newcommand{\tV}{\tens{V}}
\newcommand{\tVt}{\tV^{\T}}
\newcommand{\thV}{\wh{\tV}}
\newcommand{\thVt}{\thV^{\T}}
\newcommand{\tU}{\tens{U}}
\newcommand{\tUt}{\tU^{\T}}
\newcommand{\thU}{\wh{\tU}}
\newcommand{\thUt}{\thU^{\T}}
\newcommand{\tX}{\tens{X}}
\newcommand{\thX}{\wh{\tX}}
\newcommand{\tY}{\tens{Y}}
\newcommand{\tYt}{\tY^{\T}}
\newcommand{\thY}{\wh{\tY}}
\newcommand{\thYt}{\thY^{\T}}

\newcommand{\qprm}{q^{\prime}}

\newcommand{\tS}{\tens{S}}
\newcommand{\tSt}{\tS^{\T}}
\newcommand{\thS}{\wh{\tS}}
\newcommand{\thSt}{\thS^{\T}}
\newcommand{\teJ}{\tens{J}}
\newcommand{\theJ}{\wh{\teJ}}
\newcommand{\tI}{\tens{I}}
\newcommand{\thI}{\wh{\tI}}
\newcommand{\tC}{\tens{C}}
\newcommand{\tCt}{\tC^{\T}}
\newcommand{\thC}{\wh{\tC}}
\newcommand{\thCt}{\thC^{\T}}
\newcommand{\tG}{\tens{G}}
\newcommand{\tGt}{\tG^{\T}}
\newcommand{\thG}{\widehat{\tG}}

\newcommand{\tVe}{\tV_{\epsilon}}
\newcommand{\tVet}{\tVe^{\T}}
\newcommand{\thVe}{\widehat{\tV}_\epsilon}
\newcommand{\thVet}{\thVe^{\T}}
\newcommand{\tUe}{\tU_{\epsilon}}
\newcommand{\tUet}{\tUe^{\T}}
\newcommand{\thUe}{\widehat{\tU}_\epsilon}
\newcommand{\thUet}{\thUe^{\T}}
\newcommand{\tSe}{\tS_{\epsilon}}
\newcommand{\tSet}{\tSe^{\T}}
\newcommand{\thSe}{\widehat{\tS}_{\epsilon}}
\newcommand{\thSet}{\thSe^{\T}}

\newcommand{\tCe}{\tC_{\epsilon}}
\newcommand{\tCet}{\tCe^{\T}}
\newcommand{\thCe}{\widehat{\tC}_\epsilon}
\newcommand{\thCet}{\thCe^{\T}}

\newcommand{\tE}{\tens{E}}
\newcommand{\thE}{\widehat{\tE}}
\newcommand{\tQ}{\tens{Q}}
\newcommand{\tQt}{\tQ^{\T}}
\newcommand{\thQ}{\widehat{\tQ}}
\newcommand{\thQt}{\thQ^{\T}}
\newcommand{\hsigma}{\ensuremath{\hat{\sigma}}}
\newcommand{\tZ}{\tens{Z}}
\newcommand{\tZt}{\tZ^{\T}}
\newcommand{\thZ}{\widehat{\tZ}}
\newcommand{\thZt}{\thZ^{\T}}

\newcommand{\tW}{\tens{W}}
\newcommand{\tWt}{\tW^{\T}}
\newcommand{\thW}{\widehat{\tW}}
\newcommand{\thWt}{\thW^{\T}}

\newcommand{\tD}{\tens{D}}
\newcommand{\thD}{\widehat{\tD}}

\newcommand{\tP}{\tens{P}}
\newcommand{\tPt}{\tP^{\T}}
\newcommand{\thP}{\widehat{\tP}}
\newcommand{\thPt}{\thP^{\T}}

\newcommand{\tH}{\tens{H}}
\newcommand{\tHt}{\tH^{\T}}
\newcommand{\thH}{\widehat{\tH}}
\newcommand{\thHt}{\thH^{\T}}

\newcommand{\upbb}[2]{\ensuremath{{#1}^{(#2)}}}

\newcommand{\tWh}{\upbb{\tW}{h}}
\newcommand{\thWh}{\widehat{\upbb{\tW}{h}}}
\newcommand{\thWht}{(\thWh)^{\T}}

\newcommand{\tCh}{\upbb{\tC}{h}}
\newcommand{\tCht}{(\tCh)^{\T}}
\newcommand{\thCh}{\widehat{\upbb{\tC}{h}}}
\newcommand{\thCht}{(\thCh)^{\T}}

\newcommand{\tWl}{\upbb{\tW, \ell}}
\newcommand{\thWl}{\widehat{\upbb{\tW}{\ell}}}
\newcommand{\thWlt}{(\thWh)^{\T}}

\newcommand{\tCl}{\upbb{\tC}{\ell}}
\newcommand{\tClt}{(\tCh)^{\T}}
\newcommand{\thCl}{\widehat{\upbb{\tC}{\ell}}}
\newcommand{\thClt}{(\thCh)^{\T}}

\newcommand{\qb}{Q_{\tB}}
\newcommand{\qv}{Q_{\tV}}
\newcommand{\pb}{P_{\tB}}
\newcommand{\tpb}{\tilde{P}_{\tB}}
\newcommand{\pv}{P_{\tV}}
\newcommand{\rbb}{R_{\matb}}
\newcommand{\rvv}{R_{\matr}}

\newcommand{\mat}[1]{\mathbf{#1}}
\newcommand{\matA}{\mat{A}}
\newcommand{\matX}{\mat{X}}
\newcommand{\matR}{\mat{R}}
\newcommand{\matRt}{\matR^\T}
\newcommand{\matQ}{\mat{Q}}
\newcommand{\matE}{\mat{E}}
\newcommand{\matQt}{\matQ^{\T}}

\newcommand{\matQip}{\ensuremath{{\matQ_i}_\perp}}
\newcommand{\matQipt}{\ensuremath{\matQip^\T}}

\newcommand{\matP}{\mat{P}}
\newcommand{\matPt}{\matP^{\T}}
\newcommand{\mattP}{\tilde{\mat{P}}}
\newcommand{\mattPt}{\mattP^{\T}}
\newcommand{\matXt}{\matX^{\T}}
\newcommand{\matAt}{\matA^{\T}}
\newcommand{\bmatA}{\bar{\matA}}
\newcommand{\bmatAt}{\bmatA^{\T}}
\newcommand{\mata}{\mat{a}}
\newcommand{\matV}{\mat{V}}
\newcommand{\matVt}{\ensuremath{\matV^{\T}}}
\newcommand{\matU}{\mat{U}}
\newcommand{\matUt}{\matU^{\T}}
\newcommand{\matS}{\mat{S}}
\newcommand{\matSt}{\ensuremath{\matS^{\T}}}
\newcommand{\matT}{\mat{T}}
\newcommand{\matTt}{\ensuremath{\matT^{\T}}}
\newcommand{\mats}{\mat{s}}
\newcommand{\hmats}{\hat{\mats}}
\newcommand{\matr}{\mat{r}}
\newcommand{\maty}{\mat{y}}
\newcommand{\amatr}{\vec{\matr}}
\newcommand{\matC}{\mat{C}}
\newcommand{\matCt}{\matC^{\T}}

\newcommand{\matAr}{\mat{A}_r}

\newcommand{\matD}{\mat{D}}
\newcommand{\matd}{\mat{d}}
\newcommand{\matds}[1]{\matd^{(#1)}}
\newcommand{\ddi}{\matds{i}}
\newcommand{\ddim}{\ensuremath{D^{\xx}}}

\newcommand{\matG}{\mat{G}}
\newcommand{\matGt}{\ensuremath{\matG^{\T}}}

\newcommand{\matZ}{\mat{Z}}
\newcommand{\matZt}{\matZ^{\T}}
\newcommand{\matY}{\mat{Y}}
\newcommand{\matYt}{\matY^{\T}}

\newcommand{\matH}{\mat{H}}
\newcommand{\matHt}{\matH^\T}

\newcommand{\matLambda}{\boldsymbol{\Lambda}}
\newcommand{\matSigma}{\boldsymbol{\Sigma}}
\newcommand{\matOmega}{\boldsymbol{\Omega}}
\newcommand{\matOmegat}{\matOmega^\T}
\newcommand{\matPsi}{\boldsymbol{\Psi}}
\newcommand{\matPsit}{\matPsi^\T}

\newcommand{\matpi}{\boldsymbol{\pi}}
\newcommand{\matpit}{\matpi^\T}

\newcommand{\matK}{\mat{K}}
\newcommand{\matM}{\mathbf{M}}
\newcommand{\matMi}{\matM^{-1}}
\newcommand{\matMt}{\matM^{\T}}
\newcommand{\matW}{\mathbf{W}}
\newcommand{\matWt}{\matW^{\T}}
\newcommand{\T}{{\scriptstyle{\ensuremath{{\mat{T}}}}}}
\newcommand{\CT}{{\scriptstyle{\ensuremath{{\mat{H}}}}}}
\newcommand{\bell}{{(\ell)}}

\renewcommand{\u}{\mat{u}}
\newcommand{\ut}{\ensuremath{\u^{\T}}}
\renewcommand{\v}{\mat{v}}
\newcommand{\rrho}{\ensuremath{\bs{\rho}}}
\newcommand{\eeta}{\ensuremath{\bs{\eta}}}
\newcommand{\pphi}{\ensuremath{\bs{\varphi}}}
\newcommand{\rb}{\mat{r}}
\newcommand{\w}{\mat{w}}
\renewcommand{\b}{\mat{b}}
\newcommand{\bt}{\ensuremath{\b^{\T}}}
\newcommand{\x}{\mat{x}}
\newcommand{\h}{\mat{h}}
\newcommand{\g}{\mat{g}}
\newcommand{\e}{\mat{e}}
\newcommand{\xt}{\x^{\T}}
\newcommand{\gt}{\g^{\T}}
\newcommand{\vt}{\v^{\T}}
\newcommand{\q}{\mat{q}}
\newcommand{\qt}{\q^{\T}}
\newcommand{\p}{\mat{p}}
\newcommand{\y}{\mat{y}}
\newcommand{\bmaty}{\bar{\y}}
\newcommand{\bmatx}{\bar{\x}}
\newcommand{\z}{\mat{z}}
\newcommand{\matB}{\mat{B}}
\newcommand{\matI}{\mat{I}}
\newcommand{\matF}{\mat{F}}
\newcommand{\matFt}{\matF^{\T}}
\newcommand{\matBt}{\matB^{\T}}
\newcommand{\si}{\ensuremath{^{(i)}}}
\newcommand{\vpsi}{\boldsymbol{\psi}}
\newcommand{\vtau}{\boldsymbol{\tau}}
\newcommand{\vepsilon}{\boldsymbol{\epsilon}}
\newcommand{\vmu}{\boldsymbol{\mu}}
\newcommand{\vom}{\boldsymbol{\omega}}
\newcommand{\matGamma}{\boldsymbol{\Gamma}}
\newcommand{\TT}{\mathcal{T}}
\newcommand{\DD}{\mathcal{D}}
\newcommand{\HH}{\mathcal{H}}
\newcommand{\FF}{\mathcal{F}}
\newcommand{\LL}{\mathcal{L}}
\newcommand{\LT}{{\LL^2}}

\newcommand{\ovec}{\mathtt{sq}}
\newcommand{\reshape}{\operatorname{reshape}}

\newcommand{\matb}{\mat{b}}
\newcommand{\matbt}{\matb^{\T}}
\newcommand{\mate}{\mat{e}}
\newcommand{\matet}{\mate^{\T}}
\newcommand{\matXi}{\boldsymbol{\Xi}}

\newcommand{\bsc}{\bs{c}}
\newcommand{\bsct}{\bsc^{\T}}

\newcommand{\tlX}{\tilde{X}}
\newcommand{\tlY}{\tilde{Y}}

\newcommand{\RR}{\mathbb{R}}
\newcommand{\NN}{\mathbb{N}}
\newcommand{\CC}{\mathbb{C}}
\newcommand{\FFF}{\mathbb{F}}
\newcommand{\EE}{\mathbb{E}}
\newcommand{\ZZ}{\mathbb{Z}}
\newcommand{\sgnErr}{E_{SGN}}
\newcommand{\tsgnErr}{\tE_{SGN}}
\newcommand{\Ei}{E^{(i)}}
\newcommand{\spacex}{\boldsymbol{\mathcal{X}}}
\newcommand{\spacey}{\boldsymbol{\mathcal{Y}}}

\newcommand{\dotp}[1]{\left\langle{#1}\right\rangle}
\newcommand{\dotps}[1]{\dotp{#1}^2}
\newcommand{\FDot}[1]{\dotp{#1}_{F}}
\newcommand{\FDotS}[1]{\dotps{#1}_{F}}
\newcommand{\LTdotp}[1]{\dotp{#1}_{\LL^2}}
\newcommand{\LTdotps}[1]{\LTdotp{#1}^2}

\newcommand{\HNorm}[1]{{\left\vert\kern-0.25ex\left\vert\kern-0.25ex\left\vert{#1}\right\vert\kern-0.25ex\right\vert\kern-0.25ex\right\vert}}
\newcommand{\HNormS}[1]{\HNorm{#1}}

\newcommand{\TNorm}[1]{\|#1\|_{2}}
\newcommand{\TNormS}[1]{\TNorm{#1}^2}
\newcommand{\GNorm}[1]{\|#1\|}
\newcommand{\GNormS}[1]{\GNorm{#1}^2}

\newcommand{\GNormSd}[1]{\GNormd{#1}^2}
\newcommand{\FNorm}[1]{\|#1\|_{F}}
\newcommand{\FNormS}[1]{\FNorm{#1}^2}
\newcommand{\FNormd}[1]{\lVert#1\rVert_{F}}
\newcommand{\FNormSd}[1]{\FNormd{#1}^2}
\newcommand{\NNorm}[1]{\|#1\|_{*}}

\newcommand{\LTNorm}[1]{\|#1\|_{\LT}}
\newcommand{\LTNormS}[1]{\LTNorm{#1}^2}

\newcommand{\tNorm}[1]{\|#1\|_{top}}
\newcommand{\tNormS}[1]{\tNorm{#1}^2}

\newcommand{\war}[1]{\ensuremath{\overrightarrow{\vphantom{A}{#1}}}}
\newcommand{\bwar}[1]{\war{\bs{#1}}}

\newcommand{\trc}{\ensuremath{\operatorname{Tr}}}
\newcommand{\trace}[1]{\ensuremath{\operatorname{Tr}\left(#1 \right)}}
\newcommand{\ftr}[1]{\ensuremath{\operatorname{f-Tr}(#1)}}

\newcommand{\clr}{\noun{clr}}
\newcommand{\rclr}{\noun{rclr}}
\newcommand{\tsvdm}{\textsc{tsvdm}\xspace}
\newcommand{\tcam}{\textsc{tcam}\xspace}
\newcommand{\tca}{\textsc{tca}\xspace}
\newcommand{\tsvdmii}{\ensuremath{t}-\textsc{svdmii}\xspace}

\newcommand{\calf}{{\cal{F}}}
\newcommand{\mx}[1]{\ensuremath{\times_{#1}}}
\newcommand{\tsub}[1]{\ensuremath{\times_{#1}}}
\newcommand{\tsM}{\ensuremath{\tsub{3}\matM}}
\newcommand{\tsMinv}{\ensuremath{\tsub{3}\matM^{-1}}}
\newcommand{\ttprod}[1]{{{\star}_{{\scriptscriptstyle{\hspace{-0.1pt}#1}}}}}
\newcommand{\Mprod}{\ttprod{\matM}}
\newcommand{\mm}{\ttprod{\matM}}
\newcommand{\mmz}{\ttprod{\matM_0}}
\newcommand{\mmo}{\ttprod{\matM_1}}
\newcommand{\ff}{\ttprod{{F}}}
\newcommand{\pp}{\ttprod{{\Phi}}}
\newcommand{\muni}{\ensuremath{\Mprod{\textnormal{-unitary}}}\xspace}
\newcommand{\morth}{\ensuremath{\Mprod{\textnormal{-orthogonal}}}\xspace}
\newcommand{\pmorth}{\ensuremath{\textnormal{pseudo }\Mprod{\textnormal{-orthogonal}}}\xspace}
\newcommand{\Pmorth}{\ensuremath{\textnormal{Pseudo }\Mprod{\textnormal{-orthogonal}}}\xspace}
\newcommand{\rnk}{\ensuremath{\operatorname{rank}}}
\newcommand{\trnk}[1]{t\operatorname{-}\rnk_{#1} }
\newcommand{\Npr}{\ensuremath{N^{\prime}}}
\newcommand{\unf}[2]{\ensuremath{{#1}_{(#2)}} }
\newcommand{\supl}[1]{\ensuremath{{#1}^{[\ell]}} }
\newcommand{\fld}[2]{\ensuremath{\operatorname{fold}({#1}, #2)}}
\newcommand{\bxx}{\smash{{\setlength{\fboxsep}{0pt}\raisebox{-.5pt}{\scalebox{1.5}{\ensuremath{\xx\vphantom{A}}}}\hspace{-1pt}}}{}}

\newcommand{\bcirc}{\mathtt{bcirc}}
\newcommand{\fold}{\operatorname{fold}}

\newcommand{\XX}{\hspace*{1pt}{\textstyle\bigtimes}}

\newcommand{\OX}{{\textstyle{\bigotimes}}}
\newcommand{\Mpinv}{\mathbf{+}}
\newcommand{\mmpinv}{\ensuremath{\Mprod}\textnormal{-pseudo inverse}\xspace}

\newcommand{\fM}{\ensuremath{\mathfrak{M}} }

\algnewcommand\Input{\item[\textbf{Input:}]}%
\algnewcommand\algorithmicinput{\textbf{Input:}}
\algnewcommand\INPUT{\item[\algorithmicinput]}

\algnewcommand\Noln{\item[\hspace{28pt}]}%
\algnewcommand\algorithmicnoln{\hspace{28pt}}
\algnewcommand\NOLN{\item[\algorithmicnoln]}

\algnewcommand\algorithmicswitch{\textbf{switch}}
\algnewcommand\algorithmiccase{\textbf{case}}
\algnewcommand\algorithmicassert{\texttt{assert}}
\algnewcommand\Assert[1]{\State \algorithmicassert(#1)}

\algnewcommand\Output{\item[\textbf{Output:}]}%

\algnewcommand\Params{\item[\textbf{Parameters:}]}%
\algnewcommand\algorithmicparams{\textbf{Parameters:}}
\algnewcommand\PARAMS{\item[\algorithmicparams]}

\algnewcommand\Outputline{\item[]}%
\newcommand{\samplemode}{\emph{sample mode} }
\newcommand{\featuremode}{\emph{feature mode} }

\newcommand{\LFB}{LFB }
\newcommand{\DFB}{DFB}
\newcommand{\supp}[1]{\ensuremath{\operatorname{supp}\{#1\}}}

\newcommand{\sigmaset}{{\scalebox{1.4}{$\upsigma$}}}


\newcommand{\pbxstd}{Fig.1b\xspace}
\newcommand{\pbxtca}{Fig.1c\xspace}
\newcommand{\pbxfunnel}{Fig.1d\xspace}
\newcommand{\pbxbars}{Fig.1e\xspace}

\newcommand{\fibersctf}{Fig.fb\xspace}
\newcommand{\fiberstca}{Fig.1g\xspace}
\newcommand{\fibersfunnel}{Fig.1h\xspace}
\newcommand{\fiberstimeseries}{Fig.1i\xspace}
\newcommand{\fibersheatmap}{Fig.1j\xspace}

\newcommand{\ibdrocs}{Fig.2a\xspace}
\newcommand{\ibdloadings}{Fig.2b\xspace}
\newcommand{\ibdtca}{Fig.2c\xspace}
\newcommand{\ibdtimeseries}{Fig.2d\xspace}

\newcommand{\snydertca}{Fig.2e\xspace}
\newcommand{\snyderheatmap}{Fig.2f\xspace}

\newcommand{\Uria}[1]{\textcolor{blue}{[Uria: #1]}}
\newcommand{\ANS}[1]{\textcolor{blue}{RE: #1}}
\newcommand{\twist}{\mathtt{twist}}
\newcommand{\Haim}[1]{\textcolor{red}{[Haim: #1]}}

\definecolor{PinkC}{HTML}{ff66ff}
\definecolor{GreenC}{HTML}{33cc33}
\definecolor{BlueC}{HTML}{0099ff}
\definecolor{PurpleC}{HTML}{cc00cc}

\newcommand{\cpink}[1]{\textcolor{PinkC}{#1}}
\newcommand{\cgreen}[1]{\textcolor{GreenC}{#1}}
\newcommand{\cblue}[1]{\textcolor{BlueC}{#1}}
\newcommand{\cpurple}[1]{\textcolor{PurpleC}{#1}}

\newcommand{\matMs}{\ensuremath{\matM}^*}
\newcommand{\YY}{\ensuremath{\mathbb{Y}}}
\newcommand{\hYY}{\ensuremath{\widehat{\mathbb{Y}}}}
\newcommand{\YYr}{\ensuremath{\YY_r}}
\newcommand{\hYYr}{\ensuremath{\hYY_r}}

\newcommand{\spc}{\ensuremath{\texttt{spec}}}
\newcommand{\wc}{\ensuremath{\rightharpoonup}}

\newcommand{\doublecref}[2]{%
  \hyperref[#2]{\namecref{#1}~\labelcref*{#1}~\ref*{#2}}%
}
\newcommand{\doubleCref}[2]{%
  \hyperref[#2]{\labelcref*{#1}~\Cref{#2}}%
}
\newcommand{\ccreff}[2]{%
  \hyperref[#2]{\Cref{#2}(\labelcref*{#1})}%
}

\newif\ifuriaread 
\uriareadfalse 
\newcommand{\ureminder}[1]{%
    \Uria{%
        Set \textbackslash{}uriaread to false if you want to remove this reminder.\newline%
        #1%
    }}
\newcommand{\urif}[1]{%
    \ifuriaread%
        \ureminder{#1}%
    \else%
    \fi}

\newcommand{\alg}[1]{\ensuremath{\mathcal{#1}} }
\newcommand{\TA}{\alg{A}}
\newcommand{\TB}{\alg{B}}
\newcommand{\TM}{\alg{M}}
\newcommand{\TE}{\alg{E}}
\newcommand{\FT}{\alg{T}}
\newcommand{\FQ}{\alg{Q}}
\newcommand{\FR}{\alg{R}}

\newcommand{\ellinf}{\ell_\infty}
\newcommand{\elltwo}{\ell_2}

\newcommand{\hspc}[1]{\ensuremath{\boldsymbol{ #1}}}
\newcommand{\cH}{\ensuremath{{\alg H}}}
\newcommand{\cE}{\ensuremath{{\hspc E}}}
\newcommand{\chH}{\ensuremath{\widehat{\cH}}}
\newcommand{\chHp}{\ensuremath{\widehat{\cH}^p}}

\newcommand{\cX}{\ensuremath{{\hspc X}}}

\newcommand{\cB}{\ensuremath{{\hspc B}}}
\newcommand{\chB}{\ensuremath{\widehat{\cB}}}
\newcommand{\chBp}{\chB^p}
\newcommand{\chBpd}{(\chB^p)^*}
\newcommand{\fphi}{\FF_{\bs{\varphi}}}
\newcommand{\fphis}{{\fphi}^*}

\newcommand{\hhat}[1]{\hat{\hspace{-1pt}{#1}}}
\newcommand{\dotph}[1]{\dotp{#1}_{\cH}}
\newcommand{\dotphh}[1]{\dotp{#1}_{\chH}}
\newcommand{\dotphhp}[1]{\dotp{#1}_{\chH^p}}
\newcommand{\mpT}{\ensuremath{m \xx p \xx {\hspc T}}}
\newcommand{\mpzo}{\ensuremath{m \xx p \xx [0,1]}}
\newcommand{\cHNorm}[1]{\GNorm{#1}_{\cH}}
\newcommand{\cHNormS}[1]{\cHNorm{#1}^2}
\newcommand{\cHsNorm}[1]{\GNorm{#1}_{\cHs}}
\newcommand{\cHsNormS}[1]{\cHsNorm{#1}^2}
\newcommand{\cHsNormd}[1]{\GNormd{#1}_{\cHs}}
\newcommand{\cHsNormSd}[1]{\cHsNormd{#1}^2}

\newcommand{\dotpx}[1]{\dotp{#1}_{\cX}}
\newcommand{\cXNorm}[1]{\GNorm{#1}_{\cX}}
\newcommand{\cXNormS}[1]{\cXNorm{#1}^2}

\newcommand{\HSNorm}[1]{\GNorm{#1}_{\operatorname{HS}}}
\newcommand{\HSNormS}[1]{\HSNorm{#1}^2}
\newcommand{\dotphs}[1]{\dotp{#1}_{\operatorname{HS}}}

\newcommand{\Blt}{\ensuremath{{\hspc B}(\ell_2)}}
\newcommand{\BltNorm}[1]{\GNorm{#1}_{\Blt}}
\newcommand{\BltNormS}[1]{\BltNorm{#1}^2}

\newcommand{\ltNorm}[1]{\GNorm{#1}_{\ell_2}}
\newcommand{\ltNormS}[1]{\ltNorm{#1}^2}
\newcommand{\dotplt}[1]{\dotp{#1}_{\ell_2}}

\newcommand{\cY}{\ensuremath{{\hspc Y}}}
\newcommand{\cYd}{\cY^*}
\newcommand{\cYNorm}[1]{\GNorm{#1}_{\cY}}
\newcommand{\cYNormS}[1]{\cYNorm{#1}^2}

\newcommand{\cYp}{\cY'}
\newcommand{\cYpNorm}[1]{\GNorm{#1}_{\cYp}}
\newcommand{\cYpNormS}[1]{\cYpNorm{#1}^2}

\newcommand{\infNorm}[1]{\GNorm{#1}_{\infty}}
\newcommand{\infNormS}[1]{\infNorm{#1}^2}

\newcommand{\chHNorm}[1]{\GNorm{#1}_{\chH}}
\newcommand{\chHNormS}[1]{\chHNorm{#1}^2}

\newcommand{\chHpNorm}[1]{\GNorm{#1}_{\chHp}}
\newcommand{\chHpNormS}[1]{\chHpNorm{#1}^2}

\newcommand{\op}{\ensuremath{{\operatorname{op}}}}

\newcommand{\opNorm}[1]{\GNorm{#1}_{\op}}
\newcommand{\opNormS}[1]{\opNorm{#1}^2}
\newcommand{\opNormd}[1]{\GNormd{#1}_{\op}}
\newcommand{\opNormSd}[1]{\opNormd{#1}^2}

\newcommand{\ropNorm}[1]{\GNorm{#1}_{\FFF{-}\op}}
\newcommand{\ropNormS}[1]{\ropNorm{#1}^2}

\newcommand{\cHopNorm}[1]{\GNorm{#1}_{\cH{-}\op}}
\newcommand{\cHopNormS}[1]{\cHopNorm{#1}^2}

\newcommand{\fqopNorm}[1]{\GNorm{#1}_{\FQ{-}\op}}
\newcommand{\fqopNormS}[1]{\fqopNorm{#1}^2}

\newcommand{\cali}{\ensuremath{\hspc I}}

\newcommand{\cPo}{\ensuremath{{\hspc P}_1}}
\newcommand{\cHz}{\ensuremath{{\hspc H}_0}}
\newcommand{\cHo}{\ensuremath{{\hspc H}_1}}
\newcommand{\cHs}{\cH_{*}}
\newcommand{\cHsn}{\cH_{\divideontimes}}
\newcommand{\cHsq}[1]{\cH_{*^{#1}}}

\newcommand{\calfi}{\ensuremath{{\calf}^{-1}}}

\renewcommand{\Re}[1]{\operatorname{Re}(#1)}
\renewcommand{\Im}[1]{\operatorname{Im}(#1)}

\newcommand{\Hom}{\operatorname{Hom}}
\newcommand{\Aut}{\operatorname{Aut}}
\newcommand{\var}{\textnormal{Var}}
\newcommand{\cov}{\mathtt{cov}}
\newcommand{\cor}{\mathtt{cor}}
\newcommand{\mean}{\mathtt{mean}}

\newcommand{\fqp}{\FQ^p}
\newcommand{\chp}{\cH^p}
\newcommand{\qrnk}{\ensuremath{\FQ\operatorname{{-}rank}}}

\newcommand{\fqmpNorm}[1]{\GNorm{#1}_{\FQ^{m \xx p}}}
\newcommand{\fqmpNormS}[1]{\GNormS{#1}_{\FQ^{m \xx p}}}

\newcommand{\diag}{\operatorname{diag}}

\newcommand{\fp}{\vartriangle}
\newcommand{\Ehm}{\wh{\matE}^{(m)}}
\newcommand{\Ehp}{\wh{\matE}^{(p)}}
\newcommand{\ehpe}{\wh{e}^{(p)}}

\newcommand{\tcite}[1]{\texorpdfstring{\cite{#1}}{}}
\newcommand{\tcitepp}[2]{\texorpdfstring{\cite[#2]{#1}}{}}
\newcommand{\bsb}[1]{{\scriptscriptstyle {\hspace{0.5pt}\setlength{\fboxsep}{1pt}\fbox{\raisebox{0pt}[\dimexpr\height+0.5pt\relax][\dimexpr\depth+.51pt\relax]{#1\hspace{0.5pt}}}}}}
\newcommand{\unfs}[2]{\smash{#1_{\left( #2 \right)} }}

\newcommand{\idem}{\operatorname{idem}}
\newcommand{\iidem}{\idem^{(\operatorname{imag})}}
\newcommand{\rgen}{\Gamma_{\RR}}
\newcommand{\irgen}{\Gamma_{i\RR}}
\newcommand{\airgen}{\Gamma_{|i\RR|}}
\newcommand{\ann}{\operatorname{Ann}}
\newcommand{\lng}{\operatorname{length}}

\newcommand{\vm}[1]{\vv{\mat{#1}}}

\newcommand{\tvar}{\var_{{tmp}}}
\newcommand{\tmean}[1]{\bar{#1}}
\begin{abstract}
    A valuable feature of the tubal tensor framework is that many familiar constructions from matrix algebra carry over to tensors, including SVD and notions of rank.
    Importantly, it has been shown that for a specific family of tubal products, an Eckart-Young type theorem holds, i.e., the best low-rank approximation of a tensor under the Frobenius norm is obtained by truncating its tubal SVD.
    In this paper, we provide a complete characterization of the family of tubal products that yield an Eckart-Young type result.
    We demonstrate the practical implications of our theoretical findings by conducting experiments with video data and data-driven dynamical systems.
\end{abstract}
\maketitle

\section{Introduction}

The Eckart-Young theorem states that given a matrix $\matX$, the solution of both the following problems
\begin{equation*}
    {\min}_{\rnk(\matY) \leq r} \FNorm{\matX - \matY} , \quad \text{ and } \quad {\min}_{\rnk(\matY) \leq r} \TNorm{\matX - \matY}
\end{equation*}
is obtained by truncating the Singular Value Decomposition (SVD) of $\matX$ to its top $r$ singular values and corresponding singular vectors.
This fundamental result has placed the SVD at the core of matrix computations, with numerous applications in data science, signal processing, machine learning, and more~\cite{GOLUBVANLOAN2013,Holmes2023,Brunton2022}.
In this work, we investigate the extension of the Eckart-Young theorem to the tubal tensor framework~\cite{kilmer2008,KilmerPNAS}.

Multiway data, such as video data varying across space, time, and color channels, is naturally represented as higher-order tensors~\cite{KoldaBader2009}.
Since the Eckart-Young result holds for matrices, applying it to higher-order tensors requires flattening the tensor into a matrix, thus compromise the multiway structure of the data, which may lead to the loss of information regarding the relationships across different modes~\cite{KilmerPNAS}.

Tensor decompositions aim to generalize matrix factorizations to higher-order tensors while preserving their inherent structure~\cite{KoldaBader2009,KoldaBallard2025}.
Notable tensor decompositions include CP, HOSVD, and Tensor-Train~\cite{Hitchcock1927,Harshman1970,Tucker1966,DeLathauwerDeMoor2000,Oseledets2011}.
Each of these methods attempts to generalize the matrix SVD by approximating a tensor in Frobenius norm under various notions of tensor rank constraints. 
These decompositions lead to different notions of low-rank approximation via truncation of the decomposition, but do not provide an Eckart-Young type result, i.e., that the low rank approximation is the best one~\cite{Hillar2013,KoldaBader2009}.

The tubal framework~\cite{kilmer2008,KilmerPNAS} is based on the tubal precept: ``a tensor is a matrix of tubes''~\cite{AMDemystifying2025}, where tubes are vector space elements.
When equipped with a suitable multiplication between tubes, namely the $\mm$-product (pronounced star-M product)~\cite{KilmerPNAS}, tubes behave like scalars that can be added and multiplied.
This scalar-like structure allows for a direct carryover of many familiar definitions and properties from matrix algebra, such as matrix multiplication and factorizations, to the tubal tensor setting~\cite{Kernfeld2015}.
The $\mm$-product between tubes and tubal tensors, is determined by a choice of an invertible matrix $\matM$.
Up to certain equivalences and symmetries ~\cite{Keegan2026,dunbar2025}, each choice of $\matM$ defines a different tubal product.

It has been shown that certain choices of $\matM$, namely, the DFT matrix ~\cite{KilmerMartin11} and a nonzero multiple of a unitary matrix ~\cite{KilmerPNAS,Kernfeld2015}, guarantee an Eckart-Young type result for the resulting tSVDM (tensor SVD under the $\mm$-product), i.e., that the best low-rank approximation of a tensor is obtained by truncating its tSVDM.
There are, however, examples for which this result does not hold, thus raising the question:

%

\begingroup
    \label{q:q1}
    {\begin{flushright}
        \emph{Which tubal products yield an Eckart-Young type result for tubal tensors?}
    \end{flushright} }
\endgroup
We answer this question here by establishing the necessary and sufficient conditions on $\matM$ for Eckart-Young type optimality to hold.

\subsection{Contributions}
We fully characterize the family of $\mm$-products that yield a low-rank optimality result for the tSVDM as those defined by matrices of the form $\matM = \matD \matQ$, where $\matQ$ is a unitary matrix and $\matD$ is a real diagonal matrix.
The case of real-valued tensors is give a special attention, as the requirement of real-valuedness of the approximations imposes additional constraints possible choices of $\matM$ and the subsequent valid multiranks of tensors under such choices.
The practical implications of our theoretical results are discussed in the context of compression and feature extraction tasks.
We show that uniform scalings of unitary transforms are optimal for data compression tasks in the sense that for a given energy retention threshold, the number of rank-1 components retained in the approximation under $\matM = \matD \matQ$ is minimized when $\matD$ is a scalar multiple of the identity matrix.
Our numerical experiments with video data and data-driven dynamical systems demonstrate the advantage of having the additional degree of freedom provided by our characterization, in applications where the prioritization of certain frequencies of the data is desirable.

\section{Preliminaries and Notation}

Throughout this work, we consider constructions over a field of real $\RR$.
The field of complex numbers is denoted by $\CC$.
We use $\FFF$ to denote either $\RR$ or $\CC$.
Scalars are denoted by lowercase latin or greek letters, e.g., $\alpha,a \in \FFF$. 
The set of natural numbers is denoted by $\NN$, and for all $n \in \NN$ we define $[n] \coloneqq \{1, 2, \dots, n \}$.
Further notation and definitions are provided in the subsequent sections as we proceed.

\subsection{Linear Algebra}
Vectors and matrices are denoted by arrow accent bold lowercase, and uppercase letters, respectively, e.g., $\vm{a} \in \FFF^n, \mat{A} \in \FFF^{m \xx p}$.
The $m \xx m$ identity matrix is denoted by $\matI_m$ or simply $\matI$ when the size is clear from the context.
The transpose and conjugate-transpose of a matrix $\matA$ are denoted by $\matA^\T$ and $\matA^\CT$, respectively.
The Frobenius inner product between two matrices $\matA, \matB \in \FFF^{m \xx p}$ is defined as $\dotp{\matA, \matB}_F = \trace{\matA^\CT \matB}$. 
Correspondingly, $\FNorm{\matA} = \sqrt{\dotp{\matA, \matA}_F}$ is the Frobenius norm of $\matA$.
Let $\matA \in \FFF^{m \xx p}$ and denote by $\matA = \matU \matSigma \matV^\CT$ the (full) SVD of $\matA$, and by $\matA_q = \matU_q \matSigma_q \matV_q^\CT $ the truncation of $\matA$'s SVD to target rank $q \leq \rnk(\matA)$.

\subsection{Tensors and the Tubal Tensor Framework}
\begingroup
We write third order tensors as $\tA \in \FFF^{m \xx p \xx n}$ with entries $a_{h,j,k} \in \FFF$ indexed by triples $(h,j,k) \in [m] \xx [p] \xx [n]$.
The vectors $\tA_{:,j,k} \in \FFF^{m}$, $\tA_{h,:,k}\in \FFF^{p}$ and $\tA_{h,j,:}\in \FFF^{n}$ are $\tA$'s \textbf{column}, \textbf{row}, and \textbf{tube fibers}, respectively.
The matrices $\tA_{h,:,:} \in \FFF^{p \xx n }$, $\tA_{:,j,:} \in \FFF^{m \xx n }$ and $\tA_{:,:,k} \in \FFF^{m \xx p }$ are the \textbf{horizontal}, \textbf{lateral}, and \textbf{frontal slices} of $\tA$, respectively.

The following standard tensor definitions are taken from~\cite{KoldaBader2009,KoldaBallard2025}.
Let $\tA \in \FFF^{m \xx p \xx n}$,  the  \textbf{mode-3 unfolding} of $\tA$ is the matrix $\unfs{\tA}{3} \in \FFF^{n \xx mp}$ whose columns are the tube fibers of $\tA$. 
Let $\matX \in \FFF^{q \xx n}$, then the  \textbf{mode-3 tensor-times-matrix product} (TTM) of $\tA$ with $\matX$ is the tensor $\tA \xx_3 \matX \in \FFF^{m \xx p \xx q}$ such that $\unfs{[\tA \xx_3 \matX]}{3} = \matX \unfs{\tA}{3} \in \FFF^{q \xx mp}$, i.e., the result of multiplying each tube fiber of $\tA$ by the matrix $\matX$.
For our purposes, we extend the definition of mode-3 TTM as follows.
\begin{definition}\label{def:mode3.squeeze}
    Let $\vm{x} \in \FFF^n$ and $a \in \FFF, \matA \in \FFF^{m \xx p}$.
    Then,
    \begin{align}
        a {\xx_3} \vm{x} 
        & \in \FFF^{1 \xx 1 \xx n}, 
        && \textnormal{ where } [a {\xx_3} \vm{x}]_{ 1, 1, k} = x_k a,  
        \quad \textnormal{for all } k \in [n] \label{defeq:intro.mode_n_twisting}\\
        \matA {\xx_3} \vm{x} 
        & \in \FFF^{m \xx p \xx n}, 
        && \textnormal{ where } [\matA {\xx_3} \vm{x}]_{ :, :, k} = x_k \matA,  
        \quad \textnormal{for all } k \in [n]. \label{defeq:intro.mat.vec.twist}
    \end{align}

    The \textbf{squeeze}  operation, denoted by $\ovec$, removes all length-1 dimensions of a tensor, e.g., $\ovec(a \xx_3 \vm{x}) = a \vm{x} \in \FFF^n$ for all $a \in \FFF$ and $\vm{x} \in \FFF^n$.
\end{definition}

From \Cref{defeq:intro.mat.vec.twist}, it follows that
\begin{equation}\label{defeq:intro.tensor_twisting_decomposition}
    \tA = {\sum}_{k=1}^{n} \tA_{:,:,k} \xx_3 \vm{e}_k \quad \textnormal{  for all } \tA \in \FFF^{m \xx p \xx n}
\end{equation}
where $\vm{e}_k \in \FFF^n$ is the $k$-th standard basis vector.

\endgroup

\subsubsection{The Tubal Tensor Framework}
The tubal tensor framework views tensors as matrices~\cite{AMDemystifying2025}.
For example, an $m \xx p \xx n$ tensor, is viewed as an $m \xx p$ matrix of $n$-dimensional tube fibers.
Due to the special role of tube fibers in this framework (see \Cref{rem:tubal.scalars}), tube fibers (simply tubes, or tubal scalars), are denoted by boldsymbol lowercase letters, e.g., $\bs{a} \in \FFF^{1 \xx 1 \xx n}$.
\begin{definition}[The \texorpdfstring{$\mm$}{star M}-product \tcitepp{Kernfeld2015}{Definitions 4.1 and 4.2}.]\label{def:starM_product}
    Let $\matM \in \CC^{n \xx n}$ be invertible. 
    The \textbf{transform domain} image of $\tX \in \CC^{m \xx p \xx n}$ specified by $\matM$ is denoted by $\thX = \tX \xx_3 \matM$.
    For all $\tA \in \CC^{m \xx p \xx n}$ and $\tB \in \CC^{p \xx q \xx n}$, write
    \begin{equation*}
            \tA \mm \tB 
            = [\wh{\tA} \fp \wh{\tB}] \xx_3 \matM^{-1} \in \CC^{m \xx q \xx n}, 
            \quad \textnormal{  where } 
            [ \wh{\tA} \fp \wh{\tB} ]_{:,:,k} 
            =\wh{\tA}_{:,:,k}  \wh{\tB}_{:,:,k} \in \CC^{m \xx q}.
    \end{equation*}
    The scaling of a tensor $\tA \in \CC^{m \xx p \xx n}$ by a tube $\bs{x} \in \CC^{1 \xx 1 \xx n}$ is
    \begin{equation}\label{defeq:starM_product.scaling}
        \bs{x} \mm \tA 
        =  \thA \xx_3 \matMi\diag\left( \ovec(\hbs{x}) \right)
        =\tA \xx_3 \matMi\diag\left( \matM \ovec(\bs{x}) \right)\matM.
    \end{equation}
\end{definition}
\begin{remark}\label{rem:tubal.scalars}
    Let $\tX \in \CC^{m \xx p \xx n}$ and define $T_{\tX} : \CC^{p \xx 1 \xx n} \to \CC^{m \xx 1 \xx n}$ by $T_{\tX}(\tY) = \tX \mm \tY$ for all $\tY \in \CC^{p \xx 1 \xx n}$.
    Observe that $T_{\tX}$ is $\mm$-linear transformation, i.e., 
    $T_{\tX}(\bs{a}_1 \mm \tY_1 + \bs{a}_2 \mm \tY_2) = \bs{a}_1 \mm T_{\tX}(\tY_1) + \bs{a}_2 \mm T_{\tX}(\tY_2)$  for all $\tY_1, \tY_2 \in \CC^{p \xx 1 \xx n}$ and $\bs{a}_1, \bs{a}_2 \in \CC^{1 \xx 1 \xx n}$.
    Formally, such mappings are homomorphisms between modules over the ring of \textbf{tubal scalars}.
    See, e.g., \cite{AMDemystifying2025,Kernfeld2015} for in-depth discussions.
\end{remark}

The resulting algebraic framework arising from $\mm$-product contains the following notions listed below.
\begin{definition}[\tcitepp{Kernfeld2015}{Sections 3. and 4.}.]\label{def:matrix.mimetic.props}
    Let $\matM \in \CC^{n \xx n}$ be invertible.

    \begin{enumerate}
        \item The $\mm$-\textbf{conjugate transpose} of $\tA \in \CC^{m \xx p \xx n}$ the tensor $\tA{}^\CT \in \CC^{p \xx m \xx n}$ such that 
        \[[\thA{}^\CT]_{:,:,k} = [ \thA_{:,:,k} ]^\CT, \quad \textnormal{ for all } k \in [n].\]
        \item The $m \xx m$ \textbf{identity} tensor is $\tI_m \in \CC^{m \xx m \xx n}$ (or $\tI$ in cases where the size is clear from the context), is such that $\tI \mm \tA = \tA$, $\tB \mm \tI = \tB$ for all $\tA , \tB$ of compatible sizes.
        The \textbf{identity tube} is denoted by $\bs{e} \in \CC^{1 \xx 1 \xx n}$. 
        We have $\wh{\tI}_{:,:,k} = \matI_m$ and 
        $\wh{\bs{e}}_{1,1,k} = 1$
        for all $k \in [n]$. 
        \item Two lateral slices $\tA, \tB \in \CC^{m \xx 1 \xx n}$ are \textbf{$\mm$-orthogonal} if $\tA^\CT \mm \tB = \bs{0} \in \CC^{1 \xx 1 \xx n}.$
        \item A tensor $\tQ \in \CC^{m \xx m \xx n}$ is \textbf{$\mm$-unitary} if $\tQ^\CT \mm \tQ = \tQ \mm \tQ^\CT  = \tI_m.$
    \end{enumerate}
\end{definition}
One of the tubal framework's main features is the matrix-mimetic singular value decomposition it admits:
\begin{theorem}[tSVDM~\tcitepp{Kernfeld2015}{Theorem 5.2}]\label{thm:tSVDM}
    Let $\matM \in \CC^{n \xx n}$ be an invertible matrix.
    For any tensor $\tA \in \CC^{m \xx p \xx n}$, there exists a decomposition of the form 
    \begin{equation}\label{defeq:tSVDM}
        \tA = \tU \mm \tS \mm \tV^\CT
    \end{equation}
    where $\tU \in \CC^{m \xx m \xx n}, \tV \in \CC^{p \xx p \xx n}$ are $\mm$-unitary tensors, and $\tS \in \CC^{m \xx p \xx n}$ is a tensor whose frontal slices are all diagonal matrices.
\end{theorem}
In this work, we utilize the following definitions of tensor rank.
\begin{definition}\label{def:ranks}
    Let $\matM \in \FFF^{n \xx n}$ be an invertible linear transform, and let $\tX \in \FFF^{m \xx p \xx n}$.
    The \textbf{multirank} of $\tX$ under $\mm$ is a vector $\bs{r} = (r_1, \ldots,r_n)$ where $r_k = \rank(\wh{\tX}_{:,:,k})$ for all $k \in [n]$~{\cite[Definition 3.5]{KilmerPNAS}}.
    Write $\tX = \tU \mm \tS \mm \tV^\CT$ the tSVDM of $\tX$ (\Cref{defeq:tSVDM}), then the \textbf{t-rank} of $\tX$ under $\mm$ is the number of nonzero diagonal tubes in $\tS$~\cite[Definition 3.4]{KilmerPNAS}.
    The \textbf{implicit-rank} under $\mm$ of a tensor $\tX$ whose multirank under $\mm$ is $\bs{r} = (r_1, r_2, \ldots, r_n)$ is given by $\sum_{k=1}^n r_k$~\cite[Definition 3.6]{KilmerPNAS}.
\end{definition}
The optimality results for low-rank approximations using the tSVDM are given next.
\begin{theorem}[\tcitepp{KilmerPNAS}{Theorems 3.7 and 3.8}]\label{thm:eckart_young_tSVDM}
    Let $\matM \in \FFF^{n \xx n}$ be a nonzero multiple of a unitary matrix.
    Let $\tA \in \FFF{}^{m \xx p \xx n}$ with $\tA = \tU \mm \tS \mm \tV{}^{\CT}$.
    Then for all $r \in [\min(m,p)]$ and $\bs{r} \in [\min(m,p)]^n$ we have
    \begin{align}
        \FNorm{\tA - \tA_r} 
        &\leq  \FNorm{\tA - \tB} 
        \quad \textnormal{  for all } \tB \in \FFF^{m \xx p \xx n} \textnormal{  with t-rank } r' \leq r \nonumber \\
        \FNormS{\tA - \tA_{\bs{r}}} 
        &\leq \FNormS{\tA - \tB} 
        \quad \textnormal{  for all } \tB \in \FFF^{m \xx p \xx n} \textnormal{  such that } \rank(\wh{\tB}_{:,:,k}) \leq r_k \textnormal{  for all } k \in [n],
    \end{align}
    where $\tA_r$ and $\tA_{\bs{r}}$ are t-rank $r$ and multirank $\bs{r}$ truncations of $\tA$ defined as
    \begin{align}
        \tA_r 
        &\coloneqq  \tU_{:,1:r}  \mm  \tS_{1:r,1:r}  \mm  \tV{}^\CT_{:,1:r} 
        =  {\sum}_{j=1}^{r} \bs{s}_{j,j} \, \mm \, \tU_{:,j} \, \mm \, \tV^\CT_{:,j}  \label{defeq:trank.truncation} \\
        [\thA_{\bs{r}}]_{:,:,k}
        &\coloneqq  \thU_{:,1:r_k,k} \, \thS_{1:r_k,1:r_k,k} \, \smash{\thV}{}^\CT_{:,1:r_k,k}  \textnormal{  for all } k \in [n]. \label{defeq:multirank.truncation.hat}
    \end{align}
\end{theorem}
When $\matM$ unclear from the context, we write $\tA_r(\matM)$ and $\tA_{\bs{r}}(\matM)$ to denote the t-rank $r$ and multirank $\bs{r}$ truncations of $\tA$ under $\mm$ defined by $\matM$.

From \Cref{defeq:multirank.truncation.hat,defeq:intro.tensor_twisting_decomposition}, it follows that $\thA_{\bs{r}} = {\sum}_{k=1}^n\thU_{:,1:r_k,k} \, \thS_{1:r_k,1:r_k,k} \, \smash{\thV}{}^\CT_{:,1:r_k,k} \xx_3 \vm{e}_k$.
Then, by applying the inverse transform, we have
\begin{equation}\label{defeq:multirank.truncation}
\tA_{\bs{r}} = {\sum}_{k=1}^n\thU_{:,1:r_k,k} \, \thS_{1:r_k,1:r_k,k} \, \smash{\thV}{}^\CT_{:,1:r_k,k} \xx_3 \matMi \vm{e}_k.
\end{equation}
\Cref{def:starM_product} is stated for tensors and transforms that are both complex.
When working with real-valued data, it is common to consider real-valued transforms because they define products under which the set of real-valued tensors is closed.
However, there are examples in the literature where non-real transforms are used to define real-valued tensor products.
Most notably, the t-product~\cite{kilmer2008} uses the Discrete Fourier Transform matrix $\matF_n \in \CC^{n \xx n}$ to define a tubal product over real-valued tensors, and a corresponding tensor SVD (tSVD).
The following result from~\cite{AMDemystifying2025} specifies exactly which transforms $\matM$ can be used to a product under which the set of real tubes is closed.
\begin{lemma}[\tcitepp{AMDemystifying2025}{Lemma 2}]\label{lem:conditions_real_tubal_ring}
    Let $\matM \in \CC^{n \xx n}$ be an invertible matrix. 
    Then $\bs{a} \mm \bs{b} \in \RR^{1 \xx 1 \xx n} $ for all $\bs{a}$, $\bs{b} \in \RR^{1 \xx 1 \xx n}$, if and only if every row of $\matM$ is either real, or obtained by entry-wise complex conjugation of exactly one other row of $\matM$.
\end{lemma}

\section{Results}\label{sec:main.results}
The results in~\cite[Theorems 3.7 and 3.8]{KilmerPNAS} (\Cref{thm:eckart_young_tSVDM} here) are stated for complex tensors.
In that setting, truncations of the tSVDM are not expected to be real in general, even when the input tensor is real.
In this work, we focus on real tensors, therefore restrict our attention to tubal products under which the set of real tensors is closed, i.e., such that the conditions in \Cref{lem:conditions_real_tubal_ring} is satisfied.
Moreover, we expect that a low-rank approximations of a real tensor will be real as well.
This requirement is not trivially satisfied, and imposes further constraints on possible target multiranks of approximations.

To state our results in the most general form, we define the notion of \textit{valid multirank}.
\begin{definition}\label{def:valid.multirank}
    Let $\matM \in \CC^{n \xx n}$ be an invertible matrix satisfying the conditions in \Cref{lem:conditions_real_tubal_ring}.
    A multirank $\bs{r} \in \NN^n$ under $\mm$ is said to be \textbf{valid} if $r_j = r_{j'}$  for any $j, j' \in [n]$ such that $\matM_{j',:} = \overline{\matM_{j,:}}$.
\end{definition}
And we have the following observation.
\begin{theorem}
    Let $\matM \in \CC^{n \xx n}$ be as in \Cref{lem:conditions_real_tubal_ring}.
    Then any real tensor has a valid multirank under $\mm$.
\end{theorem}
\begin{proof}
    Let $\tX \in \RR^{m \xx p \xx n}$ and let 
    $j ,j' \in [n]$ be  such that $\matM_{j',:} = \overline{\matM_{j,:}}$.
    Then $\thX_{:,:,j} = \overline{\thX_{:,:,j'}}$, and
    thus, $\rank(\thX_{:,:,j}) = \rank(\thX_{:,:,j'})$.
\end{proof}
Therefore, our main result, stated below, is restricted to valid multiranks.
\begin{theorem}
    \label{thm:necessary_sufficient_conditions_optimality}
    For any tensor $\tX \in \RR^{m \xx p \xx n}$ and any valid multirank $\bs{r} = (r_1,\ldots,r_n)$ under $\mm$,
    \begin{equation}\label{eq:optimality.prop}
        \FNormS{\tX - \tX_{\bs{r}}} \leq \FNormS{\tX - \tY} 
        \quad \text{ for all } \tY \in \RR^{m \xx p \xx n} \text{ such that } \rank(\wh{\tY}_{:,:,k}) \leq r_k \text{ for all } k \in [n],
    \end{equation}
    holds if and only if $\matM$ can be written as
    \begin{equation}\label{eq:M.DQ.necessary_sufficient_conditions_optimality}
        \matM = \mat{D} \mat{Q},\quad \matD = \diag(d_1, \ldots, d_n) \in \RR^{n \xx n} , ~~ \matQ \in \CC^{n \xx n} , ~~
        \matQ^\CT \matQ = \mat{I} ,
    \end{equation}
    %
    where $d_1, \ldots, d_n$ are nonzero, real scalars, and the conditions of \Cref{lem:conditions_real_tubal_ring} hold for $\matM$.
\end{theorem}
We start by showing that the conditions on $\matM$ in \Cref{thm:necessary_sufficient_conditions_optimality} are sufficient for Eckart-Young optimality.
Note that this result exteds the sufficient conditions for Eckart-Young optimality in~\cite[Theorem 3.8]{KilmerPNAS} from nonzero multiple of a unitary matrix to diagonal scaling of a unitary matrix.
\begin{lemma}\label{lem:mr.trunc.optimality.sufficient}
    Let $\matM = \matD \matQ \in \CC^{n \xx n}$ be an invertible matrix satisfying \Cref{eq:M.DQ.necessary_sufficient_conditions_optimality}.
    Then, the inequality in \Cref{eq:optimality.prop} holds for any tensor $\tX \in \RR^{m \xx p \xx n}$ and a valid multirank $\bs{r} \in \NN^n$ under $\mm$.
\end{lemma}
We will use the following technical lemma.
\begin{lemma}\label{lem:DQ.frobenius_norm_expansion}
    Let $\matM = \matD \matQ$ be as in \Cref{lem:mr.trunc.optimality.sufficient}.
    Then, $\FNormS{\tX} = {\sum}_{k=1}^n d_k^{-2} \FNormS{\thX_{:,:,k}}$ for all $\tX \in \RR^{m \xx p \xx n}$.
\end{lemma}

\begin{proof}
    See \Cref{proof:DQ.frobenius_norm_expansion}.
\end{proof}

\begin{proof}[Proof of \Cref{lem:mr.trunc.optimality.sufficient}]
    Let $\tX \in \RR^{m \xx p \xx n}$.
    For any $\tY \in \RR^{m \xx p \xx n}$, we have
    \begin{align*}
        \FNormS{\tX - \tY} 
        &= {\sum}_{k=1}^n d_k^{-2} \FNormS{\thX_{:,:,k} - \thY_{:,:,k}} 
        \geq {\sum}_{k=1}^n d_k^{-2} \FNormS{\thX_{:,:,k} - [\thX_{:,:,k}]_{r_k}}
    \end{align*}
    where the first equality follows from \Cref{lem:DQ.frobenius_norm_expansion} and the inequality follows from the Eckart-Young theorem for matrices.
    From \Cref{lem:DQ.frobenius_norm_expansion} again, we have that 
    $
    {\sum}_{k=1}^n d_k^{-2} \FNormS{\thX_{:,:,k} - [\thX_{:,:,k}]_{r_k}} = \FNormS{\tX - \tX_{\bs{r}}},
    $
    and the result follows.
\end{proof}

Next, we present the proof for the necessity part of \Cref{thm:necessary_sufficient_conditions_optimality} which is the main contribution of this work.  
The proof below consider the case where $\matM \in \RR^{n \xx n}$, and the general case is established in \Cref{app:complex_case}.
\begin{proof}[Proof of necessity in \Cref{thm:necessary_sufficient_conditions_optimality} for \texorpdfstring{$\matM \in \RR^{n \xx n}$}{real transforms}]
    Let $\bs{a} \in \RR^{1 \xx 1 \xx n}$, and define $\bs{r} = \vm{e}_k$, i.e., $r_k = 1$ for some $k \in [n]$ and $r_{k'} = 0$ for all $k' \neq k$.
    The best Frobenius-norm approximation of $\bs{a}$ with multirank at most $\bs{r}$ is exactly the solution of the constrained minimization problem:
    \begin{equation}\label{eq:constrained.minimization.problem}
        {\min}_{\bs{b} \in \RR^{1 \xx 1 \xx n}} \FNormS{\bs{a} - \bs{b}} \quad \text{ subject to } \quad \rank(\hbs{b}_{:,:,k'}) \leq r_{k'} \text{ for all } k' \in [n].
    \end{equation}
    Note if a tube $\bs{b} \in \RR^{1 \xx 1 \xx n}$ has multirank at most $\bs{r} $ then $\bs{b}_{\bs{r}} = \bs{b}$.
    This allows us to rewrite \Cref{eq:constrained.minimization.problem} as 
    \begin{equation}\label{eq:constrained.minimization.problem2}
        {\min}_{\bs{b} \in \RR^{1 \xx 1 \xx n}} \FNormS{\bs{a} - \bs{b}} \quad \text{ subjected to } \quad \bs{b} = \bs{b}_{\bs{r}}.
    \end{equation}
    By \Cref{defeq:multirank.truncation,defeq:intro.mode_n_twisting}, $\bs{b} = \bs{b}_{\bs{r}} = \beta \xx_3 \matMi \vm{e}_k$.
    As a result, the problem in \Cref{eq:constrained.minimization.problem2} is equivalent to minimizing the function $F(\beta) = \FNormS{\bs{a} - \beta \xx_3 \matMi \vm{e}_k}$ over the scalar $\beta \in \RR$.
    By assumption, the solution of this problem is given by $\bs{b}_* = \bs{a}_{\bs{r}} = \wh{a}_k \xx_3 \matMi \vm{e}_k$.
    
    The first-order optimality condition implies that $F'(\wh{a}_k) = 0$.
    Rewriting $F(\beta)$ as
    \begin{align*}
        F(\beta)
        &= \dotp{(\hbs{a} - \beta \xx_3 \vm{e}_k) \xx_3 \matMi, (\hbs{a} - \beta \xx_3 \vm{e}_k) \xx_3 \matMi}_F 
        = (\ovec(\hbs{a}) - \beta \vm{e}_k)^\CT \matG (\ovec(\hbs{a}) - \beta \vm{e}_k),
    \end{align*}
    where $\matG = (\matM \matM^\CT)^{-1}$. 
    Differentiating with respect to $\beta$, we get $F'(\beta) = -2 \vm{e}_k^\CT \matG (\ovec(\hbs{a}) - \beta \vm{e}_k)$.
    Equating $F'(\wh{a}_k)$ to zero, 
    we obtain that ${\sum}_{k' \neq k} \wh{a}_{k'} g_{k',k} = 0$.
    Since the choices of $\bs{a} \in \RR^{1 \xx 1 \xx n}$ and $k$ were arbitrary, it follows that $g_{k',k} = 0$ for all $k' \neq k$.
    Thus, $\matG$ is diagonal, and therefore $\matM \matM^\CT$ is diagonal as well.
    Therefore, the rows of $\matM$ are mutually orthogonal, and the result follows.
\end{proof}

\Cref{thm:necessary_sufficient_conditions_optimality} is concerned with optimality of multirank approximations.
However, since truncations at a given t-rank are a special case of multirank truncations, the results of \Cref{thm:necessary_sufficient_conditions_optimality} naturally extend to t-rank approximations as well.

\begin{proposition}
    Let $\matM$ be as in \Cref{lem:conditions_real_tubal_ring}.
    Suppose that  \Cref{eq:optimality.prop} holds for any tensor $\tX \in \RR^{m \xx p \xx n}$ and any valid multirank $\bs{r} \in \NN^n$ under $\mm$.
    Then, for any tensor $\tA \in \RR^{m \xx p \xx n}$ and  $r \in  [\min(m,p)]$, we have 
    \begin{equation}\label{eq:t-rank.optimality.prop}
        \FNormS{\tA - \tA_{r}} \leq \FNormS{\tA - \tB} 
        \quad \text{ for all } \tB \in \RR^{m \xx p \xx n} \text{ with t-rank } r' \leq r.
    \end{equation}
\end{proposition}
\begin{proof}
    Given $r$ and $\tA$, define $\bs{r} = (r, r, \ldots, r)$.
    Note that $\bs{r}$ is valid under $\mm$. 
    By \Cref{defeq:trank.truncation,defeq:multirank.truncation}, we have that $\tA_r = \tA_{\bs{r}}$.
    Let $\tB \in \RR^{m \xx p \xx n}$ be a tensor whose t-rank is at most $r$.
    From \Cref{def:ranks}, we have $\rank(\wh{\tB}_{:,:,k}) \leq r$ for all $k \in [n]$.
    Therefore, the multirank of $\tB$ is at most $\bs{r}$.
    From the Eckart-Young optimality assumption on multirank truncations, we have $\FNorm{\tA - \tA_{r}} = \FNorm{\tA - \tA_{\bs{r}}} \leq \FNorm{\tA - \tB}$.
\end{proof}

\newcommand{\On}{O_{{\scriptstyle n}}}
\newcommand{\rgz}{\RR_{{\xx}}}
\section{Practical Implications}\label{sec:practical_implications}
The compression rates of a low-multirank truncations are determined by $\matM$.
Heuristics for choosing the `right' $\matM$ have been guided by, e.g., domain knowledge~\cite{keegan2021tensor,Mor2022}, data-dependent transformations~\cite{KilmerPNAS} and downstream computational tasks~\cite{LizKatherine2024PROJ,kong2021tensor}.

Based on existing optimality results~\cite{KilmerPNAS}, it has been common practice to restrict $\matM$ as a unitary matrix to ensure best low-rank approximation guarantees.
For example,~\cite{LizKatherine2024PROJ} suggested a Riemannian optimization based framework for learning orthogonal transformation that simultaneously searches for $(\matM, \tX) \in \On \xx \mathfrak{X}$ that minimizes least-squares regression error or low-rank approximation error, where $\On$ is the set of orthogonal $n \xx n$ matrices and $\mathfrak{X}$ is some subset of $\RR^{m \xx p \xx n}$ containing the tensors of interest.

The set of matrices identified in \Cref{thm:necessary_sufficient_conditions_optimality}, strictly contains that of nonzero scalar multiples of unitary matrices, which were shown in~\cite{KilmerPNAS} as sufficient for Eckart-Young optimality of tSVDM truncations.
Here we discuss the practical implications of this result.
For simplicity, we focus on the case of real $\matM$, but the results shown here can be easily extended to the complex case.

Let $\On$ denote the set of orthogonal $n \xx n$ matrices, and $\rgz$ the set of nonzero real numbers.
Write
\[
\rgz \On = \{ c \matW \colon c \in \rgz, \matW \in \On \} \quad \text{and} \quad \rgz^n  \On = \{ \diag(\vm{d}) \matW \colon \vm{d} \in \rgz^n, \matW \in \On \}.
\]
It is easy to verify that a real matrix $\matM$ satisfies the conditions of \Cref{thm:necessary_sufficient_conditions_optimality} if and only if $\matM \in \rgz^n  \On$, and that the conditions of \Cref{thm:eckart_young_tSVDM} are equivalent to $\matM \in \rgz \On$.

The most immediate application of Eckart-Young type results is data compression. 
In this context, the question is whether the extended search space contains matrices that are better for compression, i.e., that yield smaller approximation error for the same representation budget.

Observe the following invariance.
\begin{proposition}\label{lem:normalized.mrtrunk.eq}
    Let $\matM_0 \in \rgz \On$, and $\matM_1 = \matD \matM_0$ with $\matD = \diag(\vm{d})$ such that $\matM_1 \in \rgz^n  \On$.

    Then, for  all $\tX \in \RR^{m \xx p \xx n}$ and any multirank $\bs{r}$ under $\matM_0$ we have $\tX_{\bs{r}}(\matM_0) = \tX_{\bs{r}}(\matM_1)$. 
\end{proposition}
\begin{proof}
    \begin{align*}
        \tX_{\bs{r}}(\matM_1)
        &= {\sum}_{k=1}^n [[\tX \xx_3 \matM_1]_{:,:,k}]_{r_k} \xx_3 \matM_1^{-1} \vm{e}_k 
        = {\sum}_{k=1}^n d_k [[\tX \xx_3 \matM_0]_{:,:,k}]_{r_k} \xx_3 \matM_0^{-1} \matD^{-1} \vm{e}_k \\
        &= {\sum}_{k=1}^n [[\tX \xx_3 \matM_0]_{:,:,k}]_{r_k} \xx_3 \matM_0^{-1} \vm{e}_k 
        = \tX_{\bs{r}}(\matM_0)
    \end{align*}
\end{proof}

The storage cost of $\tX_{\bs{r}} \in \RR^{m \xx p \xx n}$ is $r(m + p)$ where $r = \sum_{k=1}^n r_k$ is the implicit rank of $\tX_{\bs{r}}$~\cite{KilmerPNAS}.
While $\tX_{\bs{r}}$ is optimal among tensors with multirank at most $\bs{r}$, it is not necessarily the best approximation of $\tX$ for a storage budget of $r(m+p)$, i.e., tensors with implicit rank at most $r$ under $\mm$.
The best per-budget approximation of a tensor is given by implicit rank truncation of the tSVDM.
\begin{definition}\label{def:implicit.rank.trunc}
    Let $\matM  \in \rgz^n  \On$, and $\tX \in \RR^{m \xx p \xx n}$ with tSVDM $\tX = \tU \mm \tS \mm \tV^\CT$.
    Write $d_k = \FNorm{\matM_{k,:}}$ for all $k \in [n]$.
    For any $r \in \NN$, define \textbf{implicit rank-$r$ truncation} of $\tX$ under $\matM$ as 
    \begin{equation}\label{eq:sum.of.rank.1.components}
        \tX_{[r]}(\matM) \coloneqq {\sum}_{h=1}^r \wh{s}_{j_h,j_h,k_h} \thU_{:, j_h,k_h} \thV_{:, j_h,k_h}^\CT \xx_3 \matM^{-1} \vm{e}_{k_h},
    \end{equation}
    where the mapping $h \mapsto(j_h,k_h)$ is a bijection such that if $h<h'$ then $d_{k_h}^{-2} \wh{s}_{j_h,j_h,k_h}^2 \geq d_{k_{h'}}^{-2} \wh{s}_{j_{h'},j_{h'},k_{h'}}^2$.
    When $\matM$ is clear from context we write $\tX_{[r]}$.
\end{definition}

\begin{theorem}\label{prop:best.implicit.rank.approx}
    Let $\matM  \in \rgz^n  \On$  and $\tX \in \RR^{m \xx p \xx n}$ be a tensor.
    For any $r $ and any tensor $\tY$ with implicit rank at most $r$, we have $\FNormS{\tX - \tX_{[r]}} \leq \FNormS{\tX - \tY}$.
\end{theorem}
\begin{proof}
    See \Cref{proof:best.implicit.rank.approx}.
\end{proof}
\begin{remark}\label{rem:tSVDMII.is.implicit.rank.trunc}
    Given $\matM = c \matW \in \rgz \On$ and a retenion threshold $\gamma \in (0,1]$, 
    the t-SVDMII~\cite[Algorithm 3]{KilmerPNAS} returns the minimal multirank $\bs{\rho}$ such that $\FNormS{\tX_{\bs{\rho}}} \geq \gamma \FNormS{\tX}$ for any input tensor $\tX$.
    The t-SVDMII operates on a tensor $\tX = \tU \mm \tS \mm \tV^\CT$ by (1) sorting the values $\{|\wh{s}_{j,j,k}|^2 \}$ in decreasing order, (2) finding $r = \min \{ r' \colon \sum_{h=1}^{r'} |\wh{s}_{j_h,j_h,k_h}|^2 \geq \gamma \FNormS{\thS} \}$ and outputs $\tX_{\bs{\rho}}$ such that $\rho_k = \max \{ j \colon |\wh{s}_{j,j,k}|^2 \geq |\wh{s}_{j_r,j_r,k_r}|^2 \}$ for all $k \in [n]$.
    Since $\matM = c \matW$, we can apply \Cref{def:implicit.rank.trunc} with $d_1 = \cdots = d_n = c$. 
    Observe that the order of singular values $\{|\wh{s}_{j,j,k}|^2 \}$ in step (1) is similar to that of $d_k^{-2} |\wh{s}_{j,j,k}|^2$ in \Cref{def:implicit.rank.trunc}. 
    Then, by \Cref{eq:sum.of.rank.1.components} we have $\tX_{\bs{\rho}} = \tX_{[r]}$.
    Thus, $\tX_{\bs{\rho}}$ is the best implicit rank-$r$ approximation of $\tX$ under $\matM$.

    If $\matM \in \rgz^n  \On$ is not a nonzero scalar multiple of an orthogonal matrix, then the order of the values $\{|\wh{s}_{j,j,k}|^2 \}$ may differ from that of $d_k^{-2} |\wh{s}_{j,j,k}|^2$. 
    In this case, the output $\tX_{\bs{\rho}}$ of the t-SVDMII may not be the best implicit rank-$r=\sum_{k=1}^n \rho_k$ approximation of $\tX$ under $\matM$.
\end{remark}
We conclude this part with the following observation.
\begin{theorem}\label{thm:imprnk.err.invariant}
    Let $\tX \in \RR^{m \xx p \xx n}$ and $r \in \NN$.
    Then, 
    \[{\inf}_{\matW \in \On} \FNormS{\tX - \tX_{[r]}(\matW)} = {\inf}_{\matM \in \rgz^n  \On} \FNormS{\tX - \tX_{[r]}(\matM)}.\]
\end{theorem}
\begin{proof}
    Since $\On$ is embedded in $\rgz^n  \On$, we have ${\inf}_{\matW \in \On} \FNormS{\tX - \tX_{[r]}(\matW)} \geq {\inf}_{\matM \in \rgz^n  \On} \FNormS{\tX - \tX_{[r]}(\matM)}$.
    Let $\{\matM_\alpha\}_{\alpha \in \NN} \subset \rgz^n  \On$ be a sequence such that $\lim_{\alpha \to \infty} \FNormS{\tX - \tX_{[r]}(\matM_\alpha)} = {\inf}_{\matM \in \rgz^n  \On} \FNormS{\tX - \tX_{[r]}(\matM)}$. 
    Then $\matM_\alpha = \diag(\vm{d}_\alpha)\matW_\alpha$ for some $\vm{d}_\alpha \in \rgz^n$ and $\matW_\alpha \in \On$ for all $\alpha$.
    Denote by $\bs{\rho}_\alpha$ the multirank of $\tX_{[r]}(\matM_\alpha)$ under $\matM_\alpha$.
    Combining \Cref{lem:normalized.mrtrunk.eq,prop:best.implicit.rank.approx}, we have $\FNormS{\tX - \tX_{[r]}(\matM_\alpha)} = \FNormS{\tX - \tX_{\bs{\rho}_\alpha}(\matW_\alpha)} \geq \FNormS{\tX - \tX_{[r]}(\matW_\alpha)}$ for all $\alpha$.
    Thus, $\inf_{\matW \in \On} \FNormS{\tX - \tX_{[r]}(\matW)} \leq \inf_{\matM \in \rgz^n  \On} \FNormS{\tX - \tX_{[r]}(\matM)}$.
\end{proof}

\Cref{thm:imprnk.err.invariant} shows that chosing $\matM$ from the larger set $\rgz^n  \On$ does not improve the per-budget approximation error. 
Therefore, we conclude that for the purpose of compression, there is no advantage in extending the search space from $\rgz \On$ to $\rgz^n  \On$.
This conclusion strengthens the theoretical justification behind methods such as ~\cite{LizKatherine2024PROJ} that use $\On$ as the search space for finding the optimal transform for low-rank approximation.

Apart from compression, there are other applications where the choice $\matM \in \rgz^n  \On$ can be beneficial.
Let $\tX \in \RR^{m \xx p \xx n}$, write $\tX = \sum_{k=1}^n \sum_{j=1}^{\rho_k} \tX{}^{(j,k)} $ where $\tX{}^{(j,k)} = \wh{s}_{j,j,k} \thU_{:, j,k} \thV{}^\CT_{:, j,k} \xx_3 \matW{}^\T \vm{e}_{k}$ is a rank-1 tensor for all $j$ and $k$, $\matW \in \On$ and $\bs{\rho} = (\rho_1, \ldots, \rho_n)$ is the multirank of $\tX$ under $\matW$.
Then, $\FNormS{\thX} = \sum_{k=1}^n \sum_{j=1}^{\rho_k} |\wh{s}_{j,j,k}|^2$.

Define $\matM = \matD \matW \in \rgz^n  \On$, and note that 
$\FNormS{\tX \xx_3 \matM} = \sum_{k=1}^n d_k^2 \sum_{j=1}^{\rho_k} |\wh{s}_{j,j,k}|^2$.
That is, the diagonal factor $\matD$ in $\matM$ re-weights the contribution of each rank-1 component in the transform domain image defined by $\matW$.
Applying the tSVDMII truncation to $\tX$ under $\matM$ results in (an) implicit rank-$r$ approximation of $\tX$, not necessarily the best one (see \Cref{rem:tSVDMII.is.implicit.rank.trunc}), but a one reflecting the re-weighting of the rank-1 components in the transform domain.
Therefore, we can use the tSVDMII, without any modifications, to amplify the contribution of certain frequencies that are relevant for the task at hand by choosing $\matM$ with appropriate diagonal factor $\matD$.
This view is reminiscent of filtering in classical signal processing~\cite{smith1997scientist}, where filters operate on signals as multiplication operators in frequency space.
Our numerical demonstrations in are focused on such applications. 

\section{Numerical Illustrations}\label{sec:numerical_illustrations}
In this section we present examples for cases where non-uniformly scaled transforms can be beneficial.
We focus on two specific applications: background subtraction in video data and tensor dynamic mode decomposition (DMD)~\cite{TDMD2025SKHT}.
Having a larger feasible set of transforms to choose from does not simplify the task of selecting an appropriate transform for a given dataset and application, which is already a challenging problem even when restricted to the family of unitary transforms~\cite{LizKatherine2024PROJ,Keegan2026,KilmerPNAS}.
In each experiment, we start with an orthogonal transform and then demonstrate how a diagonal scaling factor can be used to improve the results.
The reasoning behind the choice of the initial orthogonal transform and that of the scaling applied to it are explained in each use case.

\subsection{Background Separation in Video Data}\label{subsec:cars}
Background separation is a key technique in video analysis~\cite{brutzer2011evaluation}, serving as a preprocessing step for tasks such as automated moving object detection, tracking, and activity recognition~\cite{gutchess2001background,Sobral2014}.
A video sequence  $\tY = \{ \matY_{t} \}_{t=1}^T \subset \RR^{m  \xx 1  \xx n}$, i.e., a collection of temporally ordered frames $\matY_t$, is assumed to be a sum of background and foreground content, $\matY_t = \matY_t{\hspace{-1pt}}^{(B)} + \matY_t{\hspace{-1pt}}^{(F)}$.
A background separation model, in general, is a projection $P \colon \RR^{m  \xx 1  \xx n} \to \RR^{m \xx 1 \xx n}$ such that $P(\matY_t) \approx \matY_t{\hspace{-1pt}}^{(B)}$ for the $t$-th frame, or $P(\tY) \approx \tY{}^{(B)}$ for the video sequence $\tY$.

Using the $\mm$ product framework, we approach background separation as a feature extraction problem, in a similar spirit to the works of~\cite{Mor2022,KilmerMartin11,keegan2021tensor}.
More specifically, we work in the offline setting, where a training video $\tX$ is used to learn a background separation model $P$.
The concrete form of $P$ is a projection defined by truncation of the tSVDMII of $\tX$ to a certain threshold $\gamma$.

\newcommand{\ntfrms}{{300}}
\newcommand{\mmc}{\ttprod{\matC_n}}

In this experiment, we use highway surveillance video data from Kaggle\footnote{\url{https://www.kaggle.com/datasets/shawon10/road-traffic-video-monitoring}}.
The video is recorded by a fixed, elevated camera overlooking a motorway with vehicles moving through the frame in multiple lanes.
The original video's resultion is 768 by 1364 pixels, and it contains 1550 frames (25 FPS) in 3 color channels.
For our analysis, the color channels are mean aggregated, and the resulting 3'rd order tensor is then downsampled by a factor of 2 in all three modes resulting in $\RR^{384 \xx 775 \xx 682}$ tensor.
The first \ntfrms ~frames of the video are treated as a training set $\tX \in \RR^{384 \xx \ntfrms \xx 682}$.
For brevity, we write $m = 384$, $n = 682$, and $T = \ntfrms$.
For reference, we choose the (normalized) DCT matrix $\matC_n \in \On$, a common choice in image processing applications due to its energy concentration properties~\cite{ahmed2006discrete}.

The foreground content of our highway surveillance video consists of localized, transient features that are relatively small compared to the background.
By the concentration property above, each frame of the background is well approximated by a few low-frequency DCT components.
This leads to the working assumption that $\tX^{(B)}$ is well approximated by $\tX_{[r]}$ (\Cref{def:implicit.rank.trunc}) for some small $r$.

\begin{figure}[H]
    \centering
    \begin{minipage}{0.3\linewidth}
    \includegraphics[height=36mm]{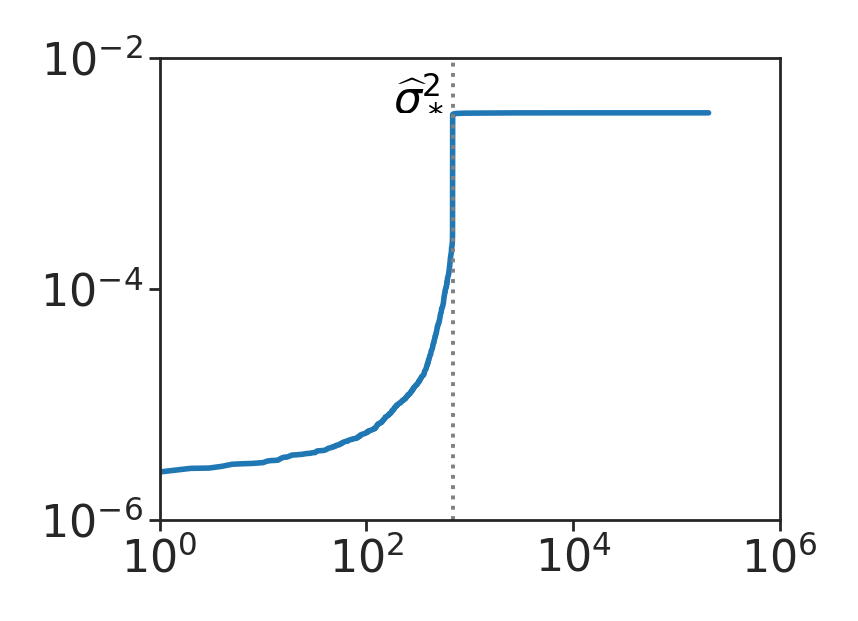}
    \end{minipage}
    \begin{minipage}{0.33\linewidth}
        \includegraphics[trim={0 6mm 0 9mm},clip,height=34mm]{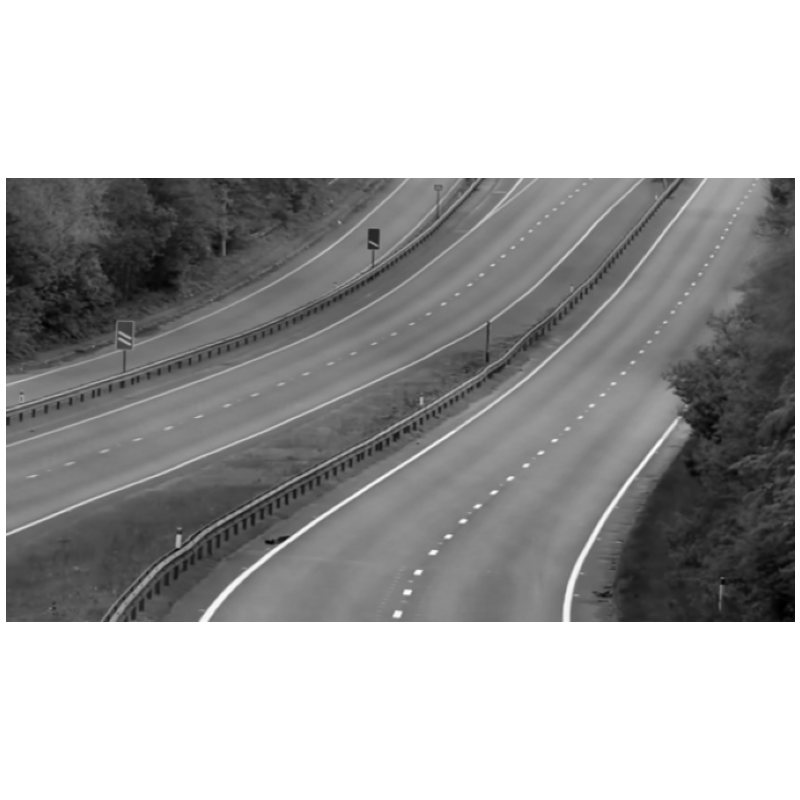}
    \end{minipage}
    \begin{minipage}{0.33\linewidth}
        \hspace{-3mm}
        \includegraphics[trim={0 6mm 0 9mm},clip,height=34mm]{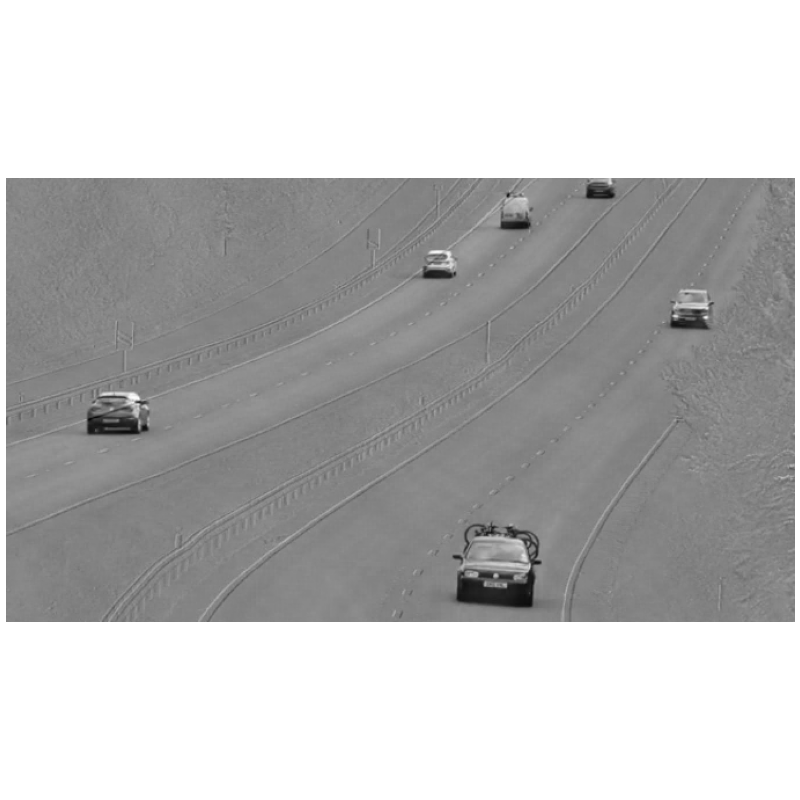}
    \end{minipage}

    \vspace{-5mm}
    \caption{\textbf{Temporal Variance Gap.} Left: y-axis shows $\wh{\sigma}^2_{\pi_i}$ for each index $i$ on the x-axis. Vertical line marks the index $i_*$ such that $\wh{\sigma}^2_{\pi_{i_*}} = \wh{\sigma}^2_*$. Images: reconstruction of visual content presented by components labaled as background (middle) and foreground (right).}
    \label{fig:var_gap}
\end{figure}

Let $\thX = \sum_{j=1}^T \sum_{k=1}^n \thX^{(j,k)}$ be the full expansion of $\thX$ as a sum of rank-1 components as in \Cref{eq:sum.of.rank.1.components}, and define the temporal variance of each component $\thX^{(j,k)}$ as
\begin{equation}\label{defeq:temporal.variance}
    \wh{\sigma}^2_{j,k}  = \frac{1}{T} {\sum}_{t=1}^T \left( \wh{v}_{t,j,k} - \frac{1}{T} {\sum}_{t'=1}^T \wh{v}_{t',j,k} \right)^2.
\end{equation}
Since the background of this video is mostly static, we can isolate it by identifying the rank-1 components with low temporal variance.
When sorted in ascending order, $\wh{\sigma}^2_{\pi_1} \leq \wh{\sigma}^2_{\pi_2} \leq \ldots \leq \wh{\sigma}^2_{\pi_{nT}}$, we see a clear gap, forming two clusters of components (\Cref{fig:var_gap}, left), which indicates a natural threshold $\wh{\sigma}^2_*$ separating the background components ($\wh{\sigma}^2_{j,k} \leq \wh{\sigma}^2_*$) from the foreground components.
\Cref{alg:var_thresh} describes a structured procedure for identifying this threshold.
\begin{algorithm}[H]
\caption{Variance threshold identification}\label{alg:var_thresh}
\begin{algorithmic}[1]
    \Input $\tX \in \RR^{m \xx T \xx n}$ (training video)
    \State Compute the rank-1 components $\{\thX^{(j,k)}\}$ from the (full) $\mmc$-tSVD of $\tX$. $\thX = \sum_{j=1}^T \sum_{k=1}^n \thX^{(j,k)}$.
    \State For each $(j,k)$, compute the temporal variance $\wh{\sigma}^2_{j,k}$ as in \Cref{defeq:temporal.variance}.
    \State Set $\pi \colon [nT] \to [T] \xx [n]$ such that $\wh{\sigma}^2_{\pi_1} \leq \cdots \leq \wh{\sigma}^2_{\pi_{nT}}$.
    \State $i_* \gets \argmax_{i=1,\ldots,T-1} \left( \wh{\sigma}^2_{\pi_{i+1}} - \wh{\sigma}^2_{\pi_i} \right)$ \Comment{index of largest gradient step}
    \Output $\wh{\sigma}^2_* = \wh{\sigma}^2_{\pi_{i_*}}$
\end{algorithmic}
\end{algorithm}
Visual inspection of the reconstructions result for each class of components confirms that their encoded visual content indeed corresponds to their assigned labels, i.e., background components correspond to the road, trees and signs (\Cref{fig:var_gap} center) and foreground components to cars (\Cref{fig:var_gap} right).

Since $\matC_n$ is orthogonal, the implicit rank-$r$ truncation of $\tX$ is equivalent to the tSVDMII truncation of $\tX$ to a suitable energy retention threshold $\gamma$ (\Cref{rem:tSVDMII.is.implicit.rank.trunc}).
We therefore seek $\gamma$ consistent with the labels just identified.

\begin{figure}[H]
    \centering
    \includegraphics[trim={0 2mm 0 5mm},clip,width=12cm]{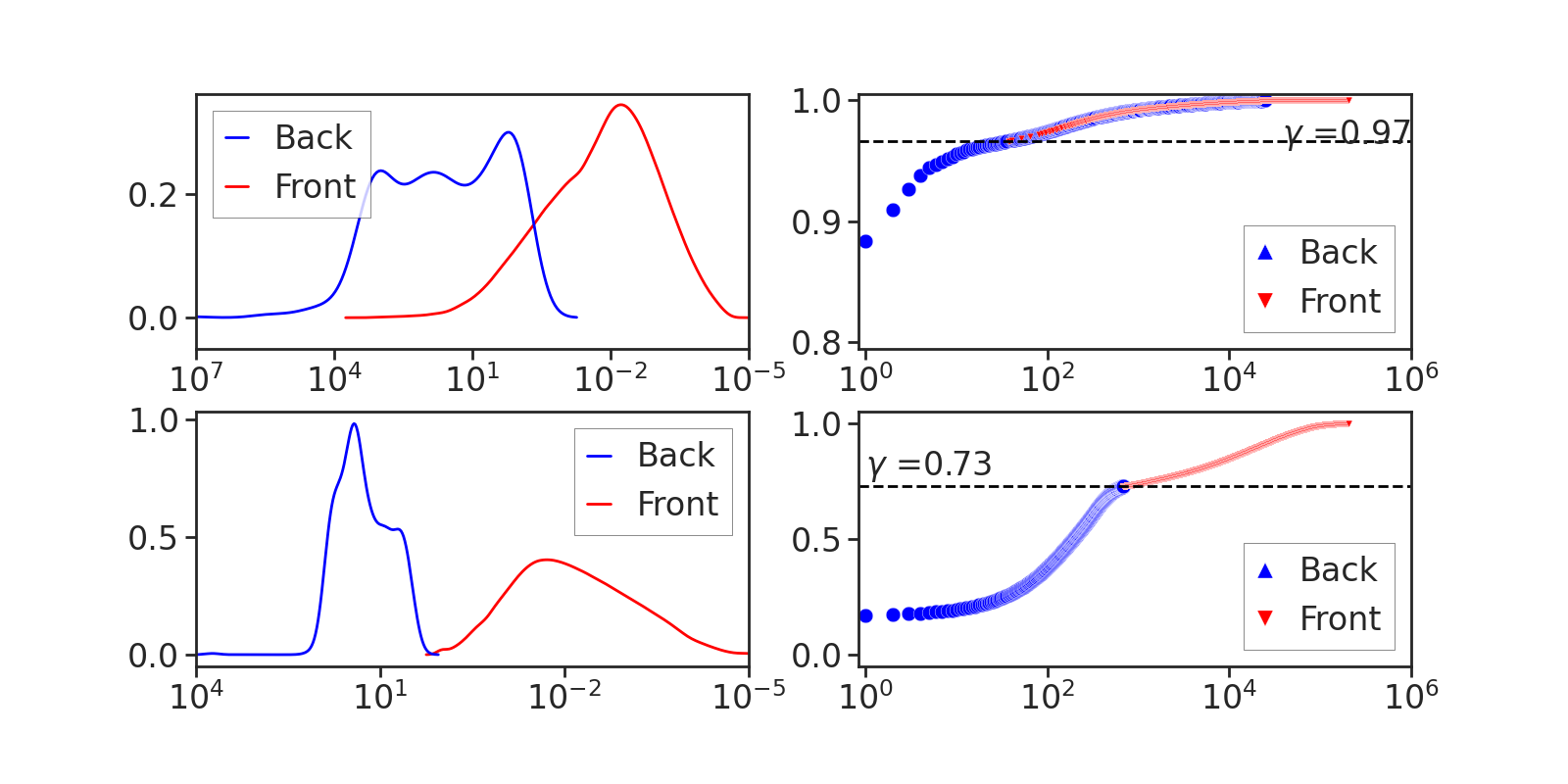}
    \caption{\textbf{Energy distribution}. Left: Kernel density estimate of $\wh{s}_{j,j,k}^2$, colored by label.
    Right: cumulative energy curves. Horizontal dashed line marks the energy retention threshold $\gamma$.
    Top: DCT, Bottom: scaled DCT.}
    \label{fig:energy_dist}
\end{figure}

We observe that background components have generally higher energy compared to foreground components (\Cref{fig:energy_dist}, top).
However, even though there is a clear separation between the two classes of components within each frequency band $k$, there is no global energy threshold that completely separates the classes.
We therefore set $\gamma$ to admit the most background-labeled energy possible without exceeding a fraction $0.01$ of the total foreground energy.
As \Cref{fig:energy_dist} (top) depicts, the resulting truncation misses a substantial portion of the background content.

A diagonal scaling $\matD = \diag(1/\wh{s}_{j_k, j_k, k})$, with $j_k = \min\{j : \wh{\sigma}^2_{j,k} \leq \wh{\sigma}^2_*\}$, normalizes each frequency band's background floor to $1$.
Since the per-frequency separation was already exact, the uniform floor lifts it to a global energy separation under $\matD\matC_n$ (\Cref{fig:energy_dist}, bottom).

\begin{figure}[H]
    \centering
    \includegraphics[width=12cm]{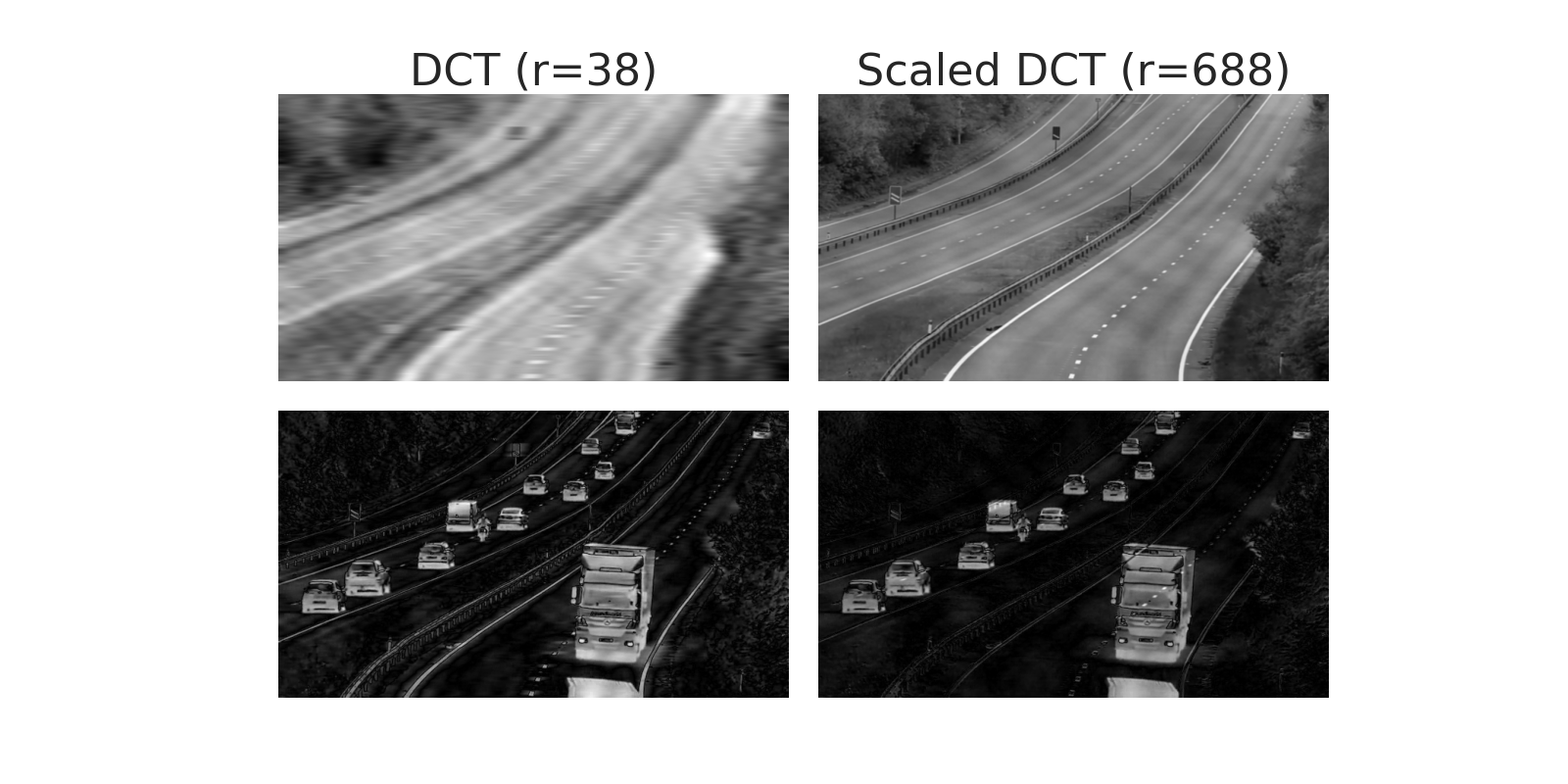}
    \caption{\textbf{Separation results}. Background (top) and foreground (bottom) estimates for the 650th frame of the original video. $r$ indicates the implicit rank of the background component.}
    \label{fig:offline_results}
\end{figure}

Next, $\tP_{\matM} = \tX_{\bs{\rho}_{\matM}} \mm \tX_{\bs{\rho}_{\matM}}^+$, where $\bs{\rho}_{\matM}$ is the multirank of the tSVDMII truncation of $\tX$ under $\matM$ to the respective energy retention threshold $\gamma$ for each choice $\matM \in \{\matC_n, \matD\matC_n\}$.
Given a test frame $\matY_t$, the background estimate is $\tP_{\matM} \mm \matY_t$, see \Cref{fig:offline_results}.

\subsection{Tensor DMD}
Consider the tensor based dynamic mode decomposition method from~\cite{TDMD2025SKHT}.
In short: given a tubal-tensor $\tX \in \RR^{m \xx p+1 \xx n}$ whose lateral slices represent $p+1$ snapshots of a dynamical system (with spatial dimension $m \xx n$), the goal is to find a linear operator $\tA_{\text{DMD}} \in \RR^{m \xx m \xx n}$ such that $\tX_{:,2:p+1,:} \approx \tA_{\text{DMD}} \mm \tX_{:,1:p,:}$.
The procedure suggested in~\cite{TDMD2025SKHT} to compute $\tA_{\text{DMD}}$ is as follows:
\begin{enumerate}
    \item Compute $\tX_{\text{train}} \coloneqq \tX_{:,1:p,:} = \tU \mm \tS \mm \tV^\CT$ the t-SVDM of $\tX_{:,1:p,:}$ under a chosen transform $\matM$.
    \item Set $\tens{K} = \tU^{\CT} \mm \tX_{:,2:p+1,:} \mm \tV \mm \tS^{+}$ (where $\tS^{+}$ is the Moore-Penrose pseudoinverse of $\tS$).
    \item Compute the Schur decomposition $\tens{K} = \tens{W} \mm \tens{T} \mm \tens{W}^{\CT}$.
    \item Set $\tZ = \tU \mm \tens{W}$
\end{enumerate}
The resulting tubal-tensor $\tZ$ contains the DMD modes of the system, and the diagonal of $\tens{T}$ contains the corresponding eigenvalues.
The DMD operator is then approximated as $\tA_{\text{DMD}} \approx \tZ \mm \tens{T} \mm \tZ^{\CT}$.
Taking low-rank truncation of $\tX_{:,1:p,:}$
in step 1 above reduces the computational and storage costs.

Similarly to~\cite{TDMD2025SKHT}, we apply the low-rank tensor DMD procedure to the cylinder flow  dataset from~\cite[Chapter 2]{Kutz2016}.
The orthogonal transformations we use as baseline are the data-driven transform, $\matZ^\T$ where $\unfs{[\tX_{\textnormal{train}}]}{3} = \matZ \matSigma \matV^\CT$~\cite{KilmerPNAS}, and the DCT matrix $\matC_n$.
The scaled versions of these transforms are defined by $\matSigma \matZ$ and $\matSigma \matC_n$. 
Under the assumption that the high energy frequencies are more relevant to the dynamics of the system, such scaling further amplifies these frequencies in comparison to the low energy frequencies.
Given implicit rank $r$, we apply the $\mm$-DMD operator $\tA_{r, \text{DMD}}$ based on truncating the t-SVDM of $\tX_{:,1:p,:}$ to $r$ components as in \Cref{eq:sum.of.rank.1.components}.

\begin{figure}[H]
    \centering
    \hspace{-10mm}\includegraphics[trim={2mm 2mm 2mm 2mm},clip,width=.87\linewidth]{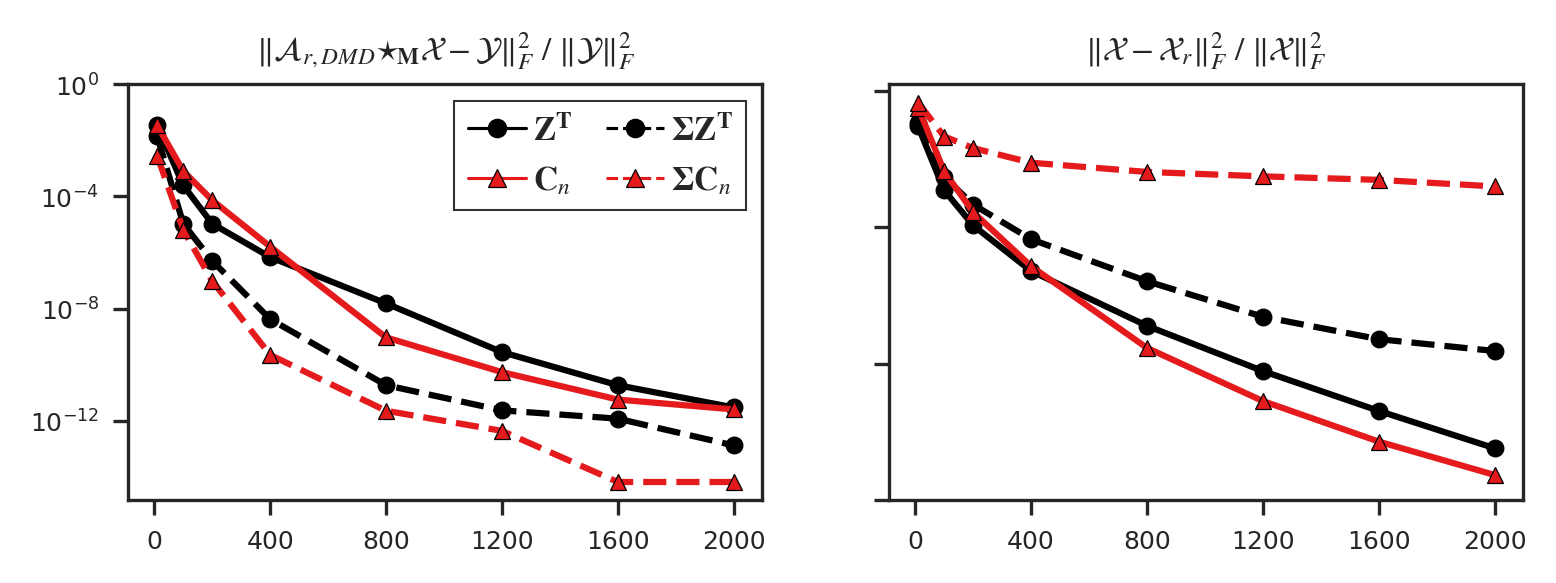}

    \includegraphics[width=.84\textwidth]{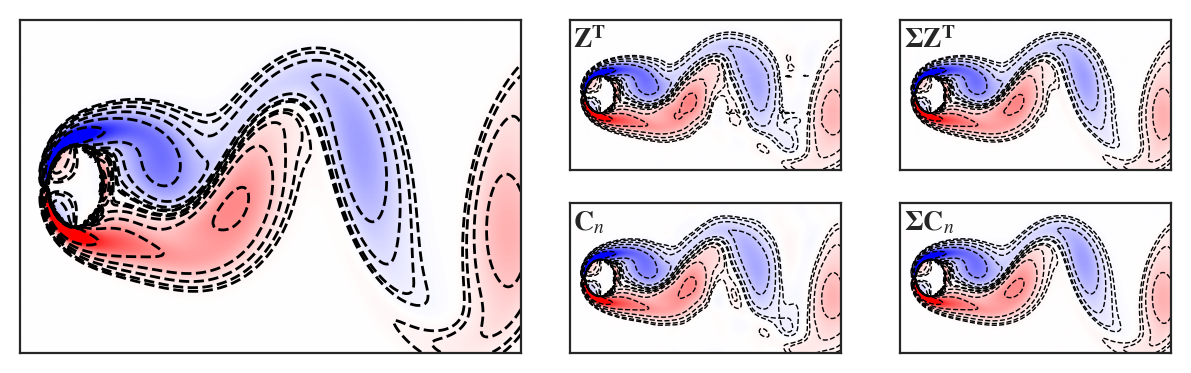}
    \caption{Top: Performance comparison of the $\mm$-DMD~\cite{TDMD2025SKHT} using orthogonal and scaled orthogonal transforms. 
    Bottom: 2d plot of flow vorticity at the last snapshots.
    Large panel: $\tX_{:,\textnormal{end},:}$. 
    Small panels show $\tA_{\text{DMD}} \mm \tX_{:,\textnormal{end},:}$.x
    Middle column: orthogonal transforms. Right column: scaled transforms.
    \label{fig1:dmdlines}}
\end{figure}

\Cref{fig1:dmdlines} shows that while the scaled transforms do not offer better compression rates of the training data, they do lead to significantly lower DMD reconstruction error compared to the unscaled transforms.
Moreover, when applying the t-SVDMII with $\gamma = 0.975$ we see
that the DMD reconstruction obtained using the orthogonal transforms contains artifacts that are not present in the DMD reconstruction obtained based on the scaled transforms.

\if11
\section{Additional Proofs}\label{sec:additional.proofs}
\subsection{\Cref{thm:necessary_sufficient_conditions_optimality} for Complex Transforms}\label{app:complex_case}


Here we prove the necessity part of
\Cref{thm:necessary_sufficient_conditions_optimality} for the general transforms $\matM \in \CC^{n \xx n}$.
Throughout this section, we write $\matG = (\matM \matM^{\CT})^{-1}$ and $g_{p,q} = \matG_{p,q}$ for any $p,q \in [n]$.
We restate the specific goal of this section for clarity.

\begin{theorem}\label{thm:necessary.conds.complex}
    Let $\matM \in \CC^{n \xx n}$ be a matrix satisfying the conditions in \Cref{lem:conditions_real_tubal_ring}.
    Suppose that the inequality in \Cref{eq:optimality.prop}: 
    \begin{equation*}
        \FNormS{\tX - \tX_{\bs{r}}} \leq \min_{\tY \in \RR^{m \xx p \xx n}} \FNormS{\tX - \tY} \quad \text{ subjected to }  \rank(\thY_{:,:,k}) \leq r_k \text{ for all } k \in [n],
    \end{equation*}
    holds for all $\tX \in \RR^{m \xx p \xx n}$ and for all valid multirank $\bs{r} = (r_1, \ldots, r_n)$.
    Then $g_{h,j} = 0 $ for any distinct indices $h,j \in [n]$.
\end{theorem}


We divide the proof of \Cref{thm:necessary.conds.complex} into three lemmas.
\Cref{lem:ey.rc.case1} shows that when $\matM_{h,:} \neq \overline{\matM_{j,:}}$, the assumptions of \Cref{thm:necessary.conds.complex} imply that $g_{h,j} = 0$.
The case where $\matM_{h,:} = \overline{\matM_{j,:}}$ is addressed by \Cref{lem:rc_case2_realpart,lem:rc_case2_impart} which show that $\Re{g_{h,j}} = 0$ and $\Im{g_{h,j}} = 0$  respectively.
We will utilize the auxiliary results stated in
\Cref{lem:transform.domain.rc.conj}.

\begin{lemma}\label{lem:transform.domain.rc.conj}
    Let $\matM \in \CC^{n \times n}$ be a matrix as in  \Cref{lem:conditions_real_tubal_ring}.
    Suppose that $k,k' \in [n]$ are such that $\matM_{k,:} = \overline{\matM_{k',:}}$, then $\thX_{:,:,k} = \overline{\vphantom{\widehat{\tX}}\thX_{:,:,k'}}$ for any $\tX \in \RR^{m \xx p \xx n}$.
    In particular, if $\Im{\matM_{k,:}} = 0 $ then $\thX_{:,:,k} = \overline{\vphantom{\widehat{\tX}}\thX_{:,:,k}} \in \RR^{m \xx p}$.
\end{lemma}
\begin{proof}
    Consider a tubal scalar $\bs{x} \in \RR^{1 \xx 1 \xx n}$, and let $k,k' \in [n]$ be such that $\matM_{k',:} = \overline{\matM_{k,:}}$.
    Write 
    \begin{align*}
        \hbs{x}_{1,1,k}
        &=\mathtt{sq}(\hbs{x})_{k} 
        = \matM_{k,:} \mathtt{sq}(\bs{x}) 
        = \overline{\matM_{k',:}} \mathtt{sq}(\bs{x}) 
        = \overline{\vm{e}_{k'}^{\CT}  \matM \mathtt{sq}(\bs{x})} 
        = \overline{\mathtt{sq}(\hbs{x})_{k'}}
        = \overline{\vphantom{\widehat{\bs{x}}}\hbs{x}_{1,1,k'}}.
    \end{align*}
    The extension to $\tX \in \RR^{m \xx p \xx n}$ follows from applying the above to each tube fiber of $\tX$.
\end{proof}



We turn to the proof of the main result of this section.
Starting with the case where $\matM_{h,:} \neq \overline{\matM_{j,:}}$.
\begin{lemma}\label{lem:ey.rc.case1}
    Let $\matM \in \CC^{n \xx n}$ be an invertible matrix satisfying \Cref{lem:conditions_real_tubal_ring}. 
    Suppose that \Cref{eq:optimality.prop} holds for any $\tX \in \RR^{m \xx p \xx n}$ and valid multirank $\bs{r}$.
    Then, $g_{h,j} = 0$ for any pair of distinct indices $h,j \in [n]$ such that $\matM_{h,:} \neq \overline{\matM_{j,:}}$.
\end{lemma}

\begin{proof}
    Define the matrix $\matP$ to be the solution of $\matP \matM = \overline{\matM}$.
    Then $\matP$ is a permutation matrix defined as
    \begin{equation}\label{eq:permutation.matrix}
        \matP = \sum_{j=1}^n \sum_{h=1}^n \xi_{j,h} \vm{e}_j  \vm{e}_h^\CT , \quad \text{ where } \xi_{j,h} = \begin{cases}
            1 & \text{ if } \matM_{h,:} = \overline{\matM_{j,:}} \\
            0 & \text{ otherwise }
        \end{cases}.
    \end{equation}

    Define 
    \begin{equation}\label{eq:highest.valid.multirank.lessthan}
        P_{\matM} (\vm{x}) = ((\matI + \matP) \vm{x} - | (\matI - \matP) \vm{x} |) / 2 ~\quad \text{ for any } \vm{x} \in \RR^n,
    \end{equation}
    where 
    $| \cdot |$ is the element-wise absolute value.
    Let $h,j \in [n]$ be distinct indices 
    and $\bs{a} \in \RR^{1 \xx 1 \xx n}$ be a tube with multirank $P_{\matM} (\vm{e}_h + \vm{e}_j)$. 
    Write $\vm{a} \coloneqq \ovec(\hbs{a})$, and note that $\vm{a}$ has the general form:
    $$
    \vm{a} = \alpha_{h} \vm{e}_{h} + \overline{\alpha_{h}} \vm{e}_{h'} + \alpha_{j} \vm{e}_{j} + \overline{\alpha_{j}} \vm{e}_{j'}
    $$
    where $\alpha_{h}, \alpha_{j} \in \CC$, $\matM_{h',:} = \overline{\matM_{h,:}}$ and $\matM_{j',:} = \overline{\matM_{j,:}}$.
    Remark that if $\Im{\matM_{h,:}} = 0$ then $h' = h$ (same for $j$).

    Let $\bs{b} \in \RR^{1 \xx 1 \xx n}$ be a tube with multirank at most $\bs{r}$, where $\bs{r} = P_{\matM} (\vm{e}_j)$ is the minimal valid multirank such that $r_j = 1$.
    Similarly to above, we have
    $
    \vm{b} \coloneqq \ovec(\hbs{b}) = \beta \vm{e}_{j} + \overline{\beta} \vm{e}_{j'}$ for some  $\beta \in \CC$.
    Express $\FNormS{\bs{a} - \bs{b}}$ as a function of the variable $\beta \in \CC$:
    \begin{align*}
        F(\beta) = \FNormS{\bs{a} - \bs{b}}
        &= \vm{a}^\CT \matG \vm{a} + (\beta \vm{e}_{j} + \overline{\beta} \vm{e}_{j'})^\CT \matG (\beta \vm{e}_{j} + \overline{\beta} \vm{e}_{j'}) 
        - \vm{a}^\CT \matG (\beta \vm{e}_{j} + \overline{\beta} \vm{e}_{j'}) 
        - (\beta \vm{e}_{j} + \overline{\beta} \vm{e}_{j'})^\CT \matG \vm{a} .
    \end{align*}
    Using identities described in~\cite{koor2023short}, we compute the Wirtinger derivative of $F$ with respect to $\beta$ as
    $$
        \frac{\partial F}{\partial \beta }
        = (\overline{\beta}-\overline{\alpha_{j}}) (g_{j,j} + g_{j',j'}) + 2({\beta} - {\alpha_{j}})  g_{j',j}
        - \overline{\alpha_{h}} (g_{h,j} + g_{j',h'}) - {\alpha_{h}} (g_{h',j} + g_{j',h}) .
    $$
    

    Supppose that the inequality in \Cref{eq:optimality.prop} holds for all $\tX \in \RR^{m \xx p \xx n}$, then $\min_\beta F(\beta) = F(\alpha_{j})$, and thus the $\frac{\partial }{\partial \beta } F\mid_{\beta = \alpha_{j}} = 0$.
    As a result, for any choice of $\alpha_{h} \in \CC$ we have
    $
        {\alpha_{h}} (g_{h',j} + g_{j',h}) = - \overline{\alpha_{h}} (g_{h,j} + g_{j',h'}) .
    $

    By taking $\alpha_{h} = 1$ we get $g_{h',j} + g_{j',h} = 0$, and $\alpha_{h} = i$ gives $g_{h,j} + g_{j',h'} = 0$.
    Note that 
    $$
        g_{j',h}
        = [\matMi_{:,j'}]^\CT \matMi_{:,h} 
        = \overline{\matMi_{:,j}}^\CT \matMi_{:,h} 
        = \overline{\matMi_{:,h}}^\CT \matMi_{:,j}
        = [\matMi_{:,h'}]^\CT \matMi_{:,j}
        = g_{h',j},
    $$
    where the second and fourth equalities follow from ~\cite[Lemma 65.]{AMDemystifying2025}. 
    Similarly, we have that $g_{h,j} = g_{j',h'}$.
    Conclude that $g_{h,j} = g_{j',h'} = g_{h',j} = g_{j',h} = 0$.
\end{proof}

\begingroup
\endgroup

The second case considers the situation where the $h,j$ are such that $\matM_{h,:} = \overline{\matM_{j,:}}$ are complex conjugates of each other.
Consider a tubal tensor $\tA \in \RR^{2 \xx 2 \xx n}$ whose frontal slices in the transform domain are zero except for the $h$-th and $j$-th slices.
\begin{align}\label{eq:Aj.Ah.case2}
    \thA_{:,:,j} &= \matA_j && \thA_{:,:,h} = \overline{\matA_j} ,
\end{align}
where $\matA_j \in \CC^{2 \xx 2}$ will be defined differently in each part of the proof.
For any $\tB \in \RR^{2 \xx 2 \xx n}$ of the same form as $\tA$, i.e., such that $\thB_{:,:,j} = \matB_j$ and  $\thB_{:,:,h} = \overline{\matB_j} $ we have,
\begin{align*}
    \FNormS{\tA - \tB}
    &= (g_{j,j} + g_{h,h})\FNormS{\matA_j - \matB_j}  
    +2\Re{ g_{h,j} \dotp{\matA_j - \matB_j,  \overline{\matA_j} - \overline{\matB_j}}_F} 
\end{align*}
For the mixed term,
\begin{align*}
    \dotp{\matA_j - \matB_j,  \overline{\matA_j} - \overline{\matB_j}}_F
    &= \FNormS{\Re{\matA_j - \matB_j}} - \FNormS{\Im{\matA_j - \matB_j}}
    + 2i \dotp{\Re{\matA_j - \matB_j},  \Im{{\matA_j} - {\matB_j}}}_F
\end{align*}
Therefore, we can write
\begin{equation}\label{eq:rc_case2.frobnorm.diff2}
    \FNormS{\tA - \tB} = S\FNormS{\matA_j - \matB_j} 
    +2\Re{ g_{h,j} \dotp{\matA_j - \matB_j,  \overline{\matA_j} - \overline{\matB_j}}_F} 
\end{equation}
where $S = g_{j,j} + g_{h,h}$.
We have the following technical lemma.
\begin{lemma}\label{lem:tech.inequalities.S.g}
    Let $\matM \in \CC^{n \xx n}$ be as in \Cref{lem:conditions_real_tubal_ring}, and let $h,j \in [n]$ be distinct indices such that $\matM_{h,:} = \overline{\matM_{j,:}}$.
    Then,

    \begin{minipage}{0.45\textwidth}
        \begin{equation}\label{eq:tech.inequality.S.g.real}
            \frac{S - \Re{ g_{h,j}}}{S-2\Re{ g_{h,j}}} > 0 
        \end{equation}
    \end{minipage}
    \begin{minipage}{0.45\textwidth}
        \begin{equation}\label{eq:tech.inequality.S.g.imag}
            \frac{S - 2 \Im{ g_{h,j}} }{S + \Im{ g_{h,j}}} > 0 
        \end{equation}
    \end{minipage}
\end{lemma}
\begin{proof}
    {For \Cref{eq:tech.inequality.S.g.real}, first note that if $\Re{ g_{h,j}} \leq 0$ then 
    $S -  2\Re{ g_{h,j}} \geq S - \Re{ g_{h,j}} $ and $S - \Re{ g_{h,j}} = S + |\Re{ g_{h,j}}| > 0$. 
    If $\Re{ g_{h,j}} \geq 0$, 
    $$S -  \Re{ g_{h,j}} \geq  S - 2 \Re{ g_{h,j}} = (\matMi_{:,h} - \matMi_{:,j})^\CT (\matMi_{:,h} - \matMi_{:,j}).$$ 
    Since $\matMi$ is invertible and $h \neq j$ we have $\matMi_{:,h} \neq \matMi_{:,j}$, and thus, $ S - 2 \Re{ g_{h,j}} > 0$.
    In both cases,  \Cref{eq:tech.inequality.S.g.real} holds

    For \Cref{eq:tech.inequality.S.g.imag}, if  $\Im{ g_{h,j}} < 0$, then $ S - 2 \Im{g_{h,j}} > S + \Im{g_{h,j}} >  S + 2\Im{g_{h,j}} $. 
    Furthermore, we have $S+ 2\Im{g_{h,j}} = (\matMi_{:,h} + i\matMi_{:,j})^\CT(\matMi_{:,h} + i\matMi_{:,j}) $.
    Since $\matMi$ is invertible and $h \neq j$, it follows that $\matMi_{:,h} \neq - i\matMi_{:,j}$  and thus $ S + 2\Im{g_{h,j}} > 0$.
    If $\Im{ g_{h,j}} > 0$, then $ S + \Im{g_{h,j}} > S - 2 \Im{g_{h,j}} > 0$, completing the proof of \Cref{eq:tech.inequality.S.g.imag}.}
\end{proof}

We split the proof and show that $\Re{g_{h,j}} = 0$ and $\Im{g_{h,j}} = 0$ separately. 

\begin{lemma}\label{lem:rc_case2_realpart}
    Let $\matM \in \CC^{n \xx n}$ be an invertible matrix satisfying \Cref{lem:conditions_real_tubal_ring}. 
    Suppose that \Cref{eq:optimality.prop} holds for any $\tX \in \RR^{m \xx p \xx n}$ and valid multirank $\bs{r}$.
    Then, $\Re{g_{h,j}} = 0$ holds for any pair of distinct indices $h,j \in [n]$ such that $\matM_{h,:} = \overline{\matM_{j,:}}$, 
\end{lemma}
\begin{proof}
    In \Cref{eq:Aj.Ah.case2}, set $\matA_j = \diag(1, \alpha_2)$ where $ \alpha_2 = i \sqrt{({S - \Re{ g_{h,j}}})/({S-2\Re{ g_{h,j}}})}$.
    Then, $\matA_j = \thA_{:,:,j}$, $\overline{\matA_j} = \diag(1, -\alpha_2) = \thA_{:,:,h}$, and $\thA_{:,:,k} = 0$ for all $k \neq h,j$.

    By \Cref{lem:tech.inequalities.S.g} and \Cref{eq:tech.inequality.S.g.real}, we have that $\alpha_2$ is well-defined and purely imaginary number with $\Im{ \alpha_2} > 0$.
    Let $\bs{r} \in \NN^n$ be a multirank vector with $r_h = r_j = 1$ and $r_k = 0$ for all $k \in [n]$ other than $h$ and $j$.
    Note that $\bs{r}$ is valid. 
    Define $\tA^{(1)}, \tA^{(2)} \in \RR^{2 \xx 2 \xx n}$ as the tensors such that 
    \begin{align*}
        \thA^{(1)}_{:,:,j} = \thA^{(1)}_{:,:,h}
        = \begin{bmatrix}
            1 & 0 \\
            0 & 0
        \end{bmatrix}, &&
        \thA^{(2)}_{:,:,j} = -\thA^{(2)}_{:,:,h}
        = \begin{bmatrix}
            0 & 0 \\
            0 & \alpha_2
        \end{bmatrix}, 
        && \text{ and } \thA^{(1)}_{:,:,k} = \thA^{(2)}_{:,:,k} = 0 \text{ for all } k \neq h,j
    \end{align*}
    From \Cref{defeq:multirank.truncation} we have 
    \begin{align}
        \thA_{\bs{r}} &= \begin{cases}
            \thA^{(1)} & \text{ if }  |\alpha_2| \leq 1  \\
            \thA^{(2)} & \text{ otherwise } 
        \end{cases}\label{eq:rc.case21.Ar.cases}
    \end{align}

    Substituting $\tB$ with $\tA^{(1)}$ and the using the concrete expression of $\tA$ in
    \Cref{eq:rc_case2.frobnorm.diff2}, we have that $\matA_j - \matB_j = \diag(0, \alpha_2)$, and thus,
    \begin{align}
        \FNormS{\tA - \tA^{(1)}} 
        = |\alpha_2|^2 (S - 2\Re{ g_{h,j}}) \label{eq:rc.case21.diffAA1}
    \end{align}
    Doing the same for $\tB = \tA^{(2)}$, we have that $\matA_j - \matB_j = \diag(1, 0)$, and therefore,
    \begin{align}
        \FNormS{\tA - \tA^{(2)}} 
        = S +2\Re{ g_{h,j} }.\label{eq:rc.case21.diffAA2}
    \end{align}
    
    Suppose that $\Re{ g_{h,j}} > 0$.
    Then $|\alpha_2| > 1$ and therefore $\thA_{\bs{r}} = \thA^{(2)}$ by \Cref{eq:rc.case21.Ar.cases}.
    Then, from \Cref{eq:optimality.prop} we have that $\FNormS{\tA - \tA^{(2)}} \leq \FNormS{\tA - \tA^{(1)}}$.
    Using \Cref{eq:rc.case21.diffAA1,eq:rc.case21.diffAA2}, we obtain
    \begin{align*}
        S +2\Re{ g_{h,j} } 
        &\leq \frac{S - \Re{ g_{h,j}}}{S-2\Re{ g_{h,j}}} (S - 2\Re{ g_{h,j}}) 
        = S - \Re{ g_{h,j}},
    \end{align*}
    therefore, $\Re{ g_{h,j}} \leq 0$ which is a contradiction.

    If $\Re{ g_{h,j}} < 0$ then  $|\alpha_2| < 1$, and by \Cref{eq:rc.case21.Ar.cases}, we have $\thA_{\bs{r}} = \thA^{(1)}$. 
    Consequently, from \Cref{eq:optimality.prop}, $\FNormS{\tA - \tA^{(1)}} \leq \FNormS{\tA - \tA^{(2)}}$.
    Applying \Cref{eq:rc.case21.diffAA1,eq:rc.case21.diffAA2} gives
    $
        S +2\Re{ g_{h,j} } \geq S - \Re{ g_{h,j}},
    $
    and hence $\Re{ g_{h,j}} \geq 0$ which is again a contradiction.
    Thus, we must have that $\Re{ g_{h,j}} = 0$.
\end{proof}
Write 
the mixed term of \Cref{eq:rc_case2.frobnorm.diff2} as
\begin{align*}
    \Re{ g_{h,j} \dotp{\matA_j - \matB_j,  \overline{\matA_j} - \overline{\matB_j}}_F} 
    &= \Re{ g_{h,j}}(\FNormS{\Re{\matA_j - \matB_j}} - \FNormS{\Im{\matA_j - \matB_j}})\\
    &- 2 \Im{ g_{h,j}} \dotp{\Re{\matA_j - \matB_j},  \Im{{\matA_j} - {\matB_j}}}_F .
\end{align*}
Then using \Cref{lem:rc_case2_realpart}, we have
As a result of \Cref{lem:rc_case2_realpart}, we have $\Re{ g_{h,j}} = 0$, and thus 
\begin{equation}\label{eq:rc_case2.frobnorm.diff3}
    \FNormS{\tA - \tB} 
    = S\FNormS{\matA_j - \matB_j}  
    - 4 \Im{ g_{h,j}} \dotp{\Re{\matA_j - \matB_j},  \Im{{\matA_j} - {\matB_j}}}_F.
\end{equation}
Proceed to the imaginary part.
\begin{lemma}\label{lem:rc_case2_impart}
    Let $\matM \in \CC^{n \xx n}$ be an invertible matrix satisfying \Cref{lem:conditions_real_tubal_ring}. 
    Suppose that \Cref{eq:optimality.prop} holds for any $\tX \in \RR^{m \xx p \xx n}$ and valid multirank $\bs{r}$.
    Then, $\Im{g_{h,j}} = 0$ holds for any pair of distinct indices $h,j \in [n]$ such that $\matM_{h,:} = \overline{\matM_{j,:}}$, 
\end{lemma}
\begin{proof}
    In \Cref{eq:Aj.Ah.case2}, set $\matA_j = \diag(a(1-i), 1+i)/\sqrt{2}$ where $a = \sqrt{ ({S - 2 \Im{g_{h,j}}})/({S {+} \Im{g_{h,j}}} }).$
    Note that, by \Cref{eq:tech.inequality.S.g.imag} $a$ is a well-defined, positive real number.
    Then $\matA_j = \thA_{:,:,j}$, $\overline{\matA_j} = \diag(a(1+i), 1-i)/\sqrt{2} = \thA_{:,:,h}$, and $\thA_{:,:,k} = 0$ for all $k \in [n] $ other than $h$ and $j$.

    Let $\bs{r} \in \NN^n$ be such that $r_h = r_j = 1$ and $r_k = 0$ for all $k$ other than $h$ and $j$.
    Define $\tA^{(1)}, \tA^{(2)} \in \RR^{2 \xx 2 \xx n}$ as the tensors such that 
    \begin{align*}
        \thA^{(1)}_{:,:,j} 
        = \begin{bmatrix}
            a \frac{1-i}{\sqrt{2}} & 0 \\
            0 & 0
        \end{bmatrix}, &&
        \thA^{(1)}_{:,:,h} 
        = \begin{bmatrix}
            a \frac{1+i}{\sqrt{2}} & 0 \\
            0 & 0
        \end{bmatrix}, &&
        \thA^{(2)}_{:,:,j} 
        = \begin{bmatrix}
            0 & 0 \\
            0 & \frac{1+i}{\sqrt{2}}
        \end{bmatrix}, &&
        \thA^{(2)}_{:,:,h} 
        = \begin{bmatrix}
            0 & 0 \\
            0 & \frac{1-i}{\sqrt{2}}
        \end{bmatrix},
    \end{align*}
    and $\thA^{(1)}_{:,:,k} = \thA^{(2)}_{:,:,k} = 0$ for all $k \in [n]$ other than $h$ and $j$.
    By \Cref{defeq:multirank.truncation},
        $\thA_{\bs{r}} = \begin{cases}
            \thA^{(1)} & \text{ if } |a| \geq 1 \\
            \thA^{(2)} & \text{ otherwise } 
        \end{cases}.$

    Plugging in $\tB = \tA^{(2)}$ in \Cref{eq:rc_case2.frobnorm.diff3}, we have $\matA_j - \matB_j = \diag(a ({1-i})/({\sqrt{2}}), 0)$, and therefore
        $
        \FNormS{\tA - \tA^{(2)}} 
        = |a|^2(S   + 2 \Im{ g_{h,j}} )$.
    Similarly, for $\tB = \tA^{(1)}$ in \Cref{eq:rc_case2.frobnorm.diff3}, we have 
    $
        \FNormS{\tA - \tA^{(1)}} 
        = S - 2 \Im{ g_{h,j}} 
    $

    If $\Im{ g_{h,j}} > 0$, then $|a| < 1$.
    By assuming \Cref{eq:optimality.prop}, we get $\FNormS{\tA - \tA^{(2)}} \leq \FNormS{\tA - \tA^{(1)}}$.
    Thus, 
    \begin{align*}
        S - 2 \Im{ g_{h,j}} 
        &\geq \FNormS{\tA - \tA^{(2)}} 
        = \frac{S - 2 \Im{ g_{h,j}}}{S +  \Im{ g_{h,j}}} (S + 2 \Im{ g_{h,j}} ) 
    \end{align*}
    This simplifies to $ 1 \geq ({S + 2 \Im{ g_{h,j}}})/({S +  \Im{ g_{h,j}}}) $, which implies that $\Im{ g_{h,j}} \leq 0$, a contradiction.

    If $\Im{ g_{h,j}} < 0$ then $|a| >1$, and then \Cref{eq:optimality.prop} implies
    that $\FNormS{\tA - \tA^{(1)}} \leq \FNormS{\tA - \tA^{(2)}}$.
    As a result,
    \begin{align*}
        S - 2 \Im{ g_{h,j}} 
        &\leq  \frac{S - 2 \Im{ g_{h,j}}}{S +  \Im{ g_{h,j}}} (S + 2 \Im{ g_{h,j}} ),
    \end{align*}
    and we get $ 1 \leq ({S + 2 \Im{ g_{h,j}}})/({S +  \Im{ g_{h,j}}}) $, thus $\Im{ g_{h,j}} \geq 0$, which is again a contradiction. 
    Therefore, we must have that $\Im{ g_{h,j}} = 0$.
\end{proof}

As discussed earlier, if $\matG$ is diagonal, then $\matM$ can be written as the product of a real-valued diagonal matrix and a unitary matrix.

\subsection{Proof of \Cref{prop:best.implicit.rank.approx}}\label{proof:best.implicit.rank.approx}
\begin{proof}
    Write $\tX = \tU \mm \tS \mm \tV^\CT$, and define $\bs{\nu} = \{ d_{k_h}^{-2} \wh{s}_{j_h,j_h,k_h} \}_{h=1}^{n \min(m,p)}$ such that $\nu_h \geq \nu_{h'}$ for all $h < h'$.
    By \Cref{defeq:multirank.truncation.hat}, 
        $\tX 
        = {\sum}_{h=1}^{n \min(m,p)} \wh{s}_{j_h,j_h,k_h} \thU_{:, j_h,k_h} \thV_{:, j_h,k_h}^\CT \xx_3 \matM^{-1} \vm{e}_{k_h}$ .
    Note that the terms in the above summation are implicit rank-1,  $\mm$-orthogonal tensors, thus, from~\cite[Proposition 38.]{AMDemystifying2025} we get $\FNormS{\tX} = {\sum}_{h=1}^{n \min(m,p)} \nu_h$.

    Let $\tY$ be a tensor with implicit rank $\rho \leq r$ under $\mm$, and denote by $\bs{\rho}$ the multirank of $\tY$.
    By \Cref{thm:necessary_sufficient_conditions_optimality} we have $\FNormS{\tX - \tX_{\bs{\rho}}} \leq \FNormS{\tX - \tY}$. 
    Observe that $\FNormS{\tX_{\bs{\rho}}} = {\sum}_{h=1}^{\rho} \nu_{\pi_h}$ for some permutation $\pi$ of $[n \min(m,p)]$.
    By the definition of $\bs{\nu}$, we have 
    ${\sum}_{h=1}^{\rho} \nu_{\pi_h} \leq {\sum}_{h=1}^{\rho} \nu_h \leq {\sum}_{h=1}^{r} \nu_h$, and 
    \[\FNormS{\tX - \tX_{[r]}} = {\sum}_{h=R+1}^{n \min(m,p)} \nu_h \leq {\sum}_{h=\rho+1}^{n \min(m,p)} \nu_{\pi_h} = \FNormS{\tX - \tX_{\bs{\rho}}} \leq \FNormS{\tX - \tY}.
    \]
\end{proof}

\subsection{Proof of \Cref{lem:DQ.frobenius_norm_expansion}}\label{proof:DQ.frobenius_norm_expansion}
\begin{proof}
    Write $\FNormS{\tX} = \dotp{\thX, \thX \xx_3 \matG}_F$
    where $\matG = (\matM \matM^\CT)^{-1}$.
    From  \Cref{defeq:intro.tensor_twisting_decomposition}, we get
    $\thX = {\sum}_{k=1}^n \thX_{:,:,k} \xx_3 \vm{e}_k$. 
    Therefore, 
    $
    \tX \xx_3 \matG = {\sum}_{k=1}^n \thX_{:,:,k} \xx_3 \vm{e}_k \xx_3 \matG = {\sum}_{k=1}^n \thX_{:,:,k} \xx_3 \matG \vm{e}_k. 
    $
    From \Cref{eq:M.DQ.necessary_sufficient_conditions_optimality}, it follows that $\matM \matM^\CT = \matD^2$. 
    Consequently, $\matG = \matD^{-2}$, and we have 
    \begin{align*}
        \FNormS{\tX} 
        &= \dotp{{\sum}_{k=1}^n \thX_{:,:,k} \xx_3 \vm{e}_k, {\sum}_{k'=1}^n \thX_{:,:,k'} \xx_3 (d_{k'}^{-2} \vm{e}_{k'})}_F 
        = {\sum}_{k=1}^n {\sum}_{k'=1}^n d_{k'}^{-2}  \dotp{\thX_{:,:,k} \xx_3 \vm{e}_k, \thX_{:,:,k'} \xx_3  \vm{e}_{k'}}_F.
    \end{align*}
    By the same reasoning as above, $\dotp{\matA \xx_3 \vm{e}_k, \matB \xx_3 \vm{e}_{k'} }_F = \dotp{\matA, \matB \xx_3 \vm{e}_{k}^\CT \vm{e}_{k'} }_F$ for all matrices $\matA, \matB$ of compatible sizes.
    In particular, $\dotp{\matA \xx_3 \vm{e}_k, \matB \xx_3 \vm{e}_{k'} }_F = 0$ for all $k' \neq k$.
    Thus, we obtain
    \begin{align*}
        \FNormS{\tX} 
        &= {\sum}_{k=1}^n d_k^{{-}2} \dotp{\thX_{:,:,k} \xx_3 \vm{e}_k, \thX_{:,:,k} \xx_3  \vm{e}_{k}}_F
        = {\sum}_{k=1}^n d_k^{{-}2} \FNormS{\thX_{:,:,k}}.
    \end{align*}
\end{proof}

\fi
\section*{Acknowledgments}
I'd like to thank my advisor, Prof. Haim Avron, for his invaluable guidance and support throughout this project.
This research was funded by the Israel Science Foundation (grant 1524/23). 
\bibliographystyle{abbrv}
\bibliography{refs}

@article{keegan2021tensor,
  title={A tensor SVD-based classification algorithm applied to fmri data},
  author={Keegan, Katherine and Vishwanath, Tanvi and Xu, Yihua},
  journal={arXiv preprint arXiv:2111.00587},
  year={2021}
}

@article{KilmerMartin11,
    title = {{Factorization strategies for third-order tensors}},
    year = {2011},
    journal = {Linear Algebra and its Applications},
    author = {Kilmer, Misha E. and Martin, Carla D.},
    number = {3},
    month = {8},
    pages = {641--658},
    volume = {435},
    publisher = {North-Holland},
    url = {https://www.sciencedirect.com/science/article/pii/S0024379510004830},
    doi = {10.1016/J.LAA.2010.09.020},
    issn = {0024-3795}
}

@article{Hillar2013,
  doi = {10.1145/2512329},
  url = {https://doi.org/10.1145/2512329},
  year = {2013},
  month = nov,
  publisher = {Association for Computing Machinery ({ACM})},
  volume = {60},
  number = {6},
  pages = {1--39},
  author = {Christopher J. Hillar and Lek-Heng Lim},
  title = {Most Tensor Problems Are {NP}-Hard},
  journal = {Journal of the {ACM}}
}

@article{Harshman1970,
author = {Harshman, Richard a},
journal = {UCLA Working Papers in Phonetics},
number = {10},
title = {{Foundations of the PARAFAC procedure: Models and conditions for an “explanatory” multimodal factor analysis}},
volume = {16},
year = {1970}
}

@article{Tucker1966,
  doi = {10.1007/bf02289464},
  url = {https://doi.org/10.1007/bf02289464},
  year = {1966},
  month = sep,
  publisher = {Springer Science and Business Media {LLC}},
  volume = {31},
  number = {3},
  pages = {279--311},
  author = {Ledyard R Tucker},
  title = {Some mathematical notes on three-mode factor analysis},
  journal = {Psychometrika}
}

@article{Oseledets2011,
author = {Oseledets, I. V.},
title = {Tensor-Train Decomposition},
journal = {SIAM Journal on Scientific Computing},
volume = {33},
number = {5},
pages = {2295-2317},
year = {2011},
doi = {10.1137/090752286},
URL = { https://doi.org/10.1137/090752286},
eprint = {https://doi.org/10.1137/090752286}
}

@book{smith1997scientist,
  title={The scientist and engineer's guide to digital signal processing},
  author={Smith, Steven W and others},
  year={1997},
  publisher={California Technical Pub. San Diego}
}

@article{KilmerPNAS,
author = {Kilmer, Misha E. and Horesh, Lior and Avron, Haim and Newman, Elizabeth},
doi = {10.1073/PNAS.2015851118},
issn = {0027-8424},
journal = {Proceedings of the National Academy of Sciences},
month = jul,
number = {28},
pages = {e2015851118},
pmid = {34234014},
publisher = {National Academy of Sciences},
title = {{Tensor-tensor algebra for optimal representation and compression of multiway data}},
url = {https://www.pnas.org/content/118/28/e2015851118 https://www.pnas.org/content/118/28/e2015851118.abstract},
volume = {118},
year = {2021}
}

@article{KoldaBader2009,
  doi = {10.1137/07070111x},
  url = {https://doi.org/10.1137/07070111x},
  year = {2009},
  month = aug,
  publisher = {Society for Industrial {\&} Applied Mathematics ({SIAM})},
  volume = {51},
  number = {3},
  pages = {455--500},
  author = {Tamara G. Kolda and Brett W. Bader},
  title = {Tensor Decompositions and Applications},
  journal = {{SIAM} Review}
}

@article{DeLathauwerDeMoor2000,
author = {De Lathauwer, Lieven and De Moor, Bart and Vandewalle, Joos},
title = {On the Best Rank-1 and Rank-(R1 ,R2 ,. . .,RN) Approximation of Higher-Order Tensors},
journal = {SIAM Journal on Matrix Analysis and Applications},
volume = {21},
number = {4},
pages = {1324-1342},
year = {2000},
doi = {10.1137/S0895479898346995},
URL = {https://doi.org/10.1137/S0895479898346995},
eprint = {https://doi.org/10.1137/S0895479898346995}
}

@article{Kernfeld2015,
  doi = {10.1016/j.laa.2015.07.021},
  url = {https://doi.org/10.1016/j.laa.2015.07.021},
  year = {2015},
  month = nov,
  publisher = {Elsevier {BV}},
  volume = {485},
  pages = {545--570},
  author = {Eric Kernfeld and Misha Kilmer and Shuchin Aeron},
  title = {Tensor{\textendash}tensor products with invertible linear transforms},
  journal = {Linear Algebra and its Applications}
}

@article{Mor2022,
  doi = {10.1371/journal.pcbi.1010212},
  url = {https://doi.org/10.1371/journal.pcbi.1010212},
  year = {2022},
  month = jul,
  publisher = {Public Library of Science ({PLoS})},
  volume = {18},
  number = {7},
  pages = {e1010212},
  author = {Uria Mor and Yotam Cohen and Rafael Vald{\'{e}}s-Mas and Denise Kviatcovsky and Eran Elinav and Haim Avron},
  editor = {Carl Herrmann},
  title = {Dimensionality reduction of longitudinal 'omics data using modern tensor factorizations},
  journal = {{PLOS} Computational Biology}
}

@misc{dunbar2025,
      title={Tensor-Tensor Products, Group Representations, and Semidefinite Programming}, 
      author={Alex Dunbar and Elizabeth Newman},
      year={2025},
      eprint={2507.12729},
      archivePrefix={arXiv},
      primaryClass={math.OC},
      url={https://arxiv.org/abs/2507.12729}, 
}

@misc{kilmer2008,
  title={A third-order generalization of the matrix svd as a product of third-order tensors},
  author={Kilmer, Misha E and Martin, Carla D and Perrone, Lisa},
  journal={Tufts University, Department of Computer Science, Tech. Rep. TR-2008-4},
  year={2008}
}

@book{Kutz2016,
  title = {Dynamic Mode Decomposition: Data-Driven Modeling of Complex Systems},
  ISBN = {9781611974508},
  url = {http://dx.doi.org/10.1137/1.9781611974508},
  DOI = {10.1137/1.9781611974508},
  publisher = {Society for Industrial and Applied Mathematics},
  author = {Kutz,  J. Nathan and Brunton,  Steven L. and Brunton,  Bingni W. and Proctor,  Joshua L.},
  year = {2016},
  month = nov 
}

@misc{AMDemystifying2025,
  doi = {10.48550/ARXIV.2506.03311},
  url = {https://arxiv.org/abs/2506.03311},
  author = {Avron,  Haim and Mor,  Uria},
  title = {Demystifying Tubal Tensor Algebra},
  publisher = {arXiv},
  year = {2025},
  copyright = {arXiv.org perpetual,  non-exclusive license}
}

@BOOK{KoldaBallard2025,
  title     = {Tensor decompositions for data science},
  author    = {Ballard, Grey and Kolda, Tamara G},
  publisher = {Cambridge University Press},
  month     =  {jun},
  year      =  {2025},
  address   = {Cambridge, England},
  ISBN = {9781009471671}
}

@book{Holmes2023,
  title = {Introduction to Scientific Computing and Data Analysis},
  ISBN = {9783031224300},
  ISSN = {2197-179X},
  url = {http://dx.doi.org/10.1007/978-3-031-22430-0},
  DOI = {10.1007/978-3-031-22430-0},
  journal = {Texts in Computational Science and Engineering},
  publisher = {Springer International Publishing},
  author = {Holmes,  Mark H.},
  year = {2023}
}

@book{GOLUBVANLOAN2013,
author = {Golub, Gene H. and Van Loan, Charles F.},
title = {Matrix Computations - 4th Edition},
publisher = {Johns Hopkins University Press},
year = {2013},
doi = {10.1137/1.9781421407944},
address = {Philadelphia, PA},
edition   = {},
URL = {https://epubs.siam.org/doi/abs/10.1137/1.9781421407944},
eprint = {https://epubs.siam.org/doi/pdf/10.1137/1.9781421407944}
}

@article{Hitchcock1927,
  title = {The Expression of a Tensor or a Polyadic as a Sum of Products},
  volume = {6},
  ISSN = {0097-1421},
  url = {http://dx.doi.org/10.1002/sapm192761164},
  DOI = {10.1002/sapm192761164},
  number = {1–4},
  journal = {Journal of Mathematics and Physics},
  publisher = {Wiley},
  author = {Hitchcock,  Frank L.},
  year = {1927},
  month = apr,
  pages = {164–189}
}

@misc{TDMD2025SKHT,
  doi = {10.48550/ARXIV.2508.10126},
  url = {https://arxiv.org/abs/2508.10126},
  author = {Saibaba,  Arvind K. and Kilmer,  Misha E. and Hall-Hooper,  Khalil and Tian,  Fan and Mize,  Alex},
  title = {A tensor-based dynamic mode decomposition based on the $\star_{\boldsymbol{M}}$-product},
  publisher = {arXiv},
  year = {2025},
  copyright = {Creative Commons Attribution 4.0 International}
}

@book{Brunton2022,
  title = {Data-Driven Science and Engineering: Machine Learning,  Dynamical Systems,  and Control},
  ISBN = {9781009098489},
  url = {http://dx.doi.org/10.1017/9781009089517},
  DOI = {10.1017/9781009089517},
  publisher = {Cambridge University Press},
  author = {Brunton,  Steven L. and Kutz,  J. Nathan},
  year = {2022},
  month = may 
}

@article{ahmed2006discrete,
  title={Discrete cosine transform},
  author={Ahmed, Nasir and Natarajan, T\_ and Rao, Kamisetty R},
  journal={IEEE transactions on Computers},
  volume={100},
  number={1},
  pages={90--93},
  year={2006},
  publisher={IEEE}
}

@article{Sobral2014,
  title = {A comprehensive review of background subtraction algorithms evaluated with synthetic and real videos},
  volume = {122},
  ISSN = {1077-3142},
  url = {http://dx.doi.org/10.1016/j.cviu.2013.12.005},
  DOI = {10.1016/j.cviu.2013.12.005},
  journal = {Computer Vision and Image Understanding},
  publisher = {Elsevier BV},
  author = {Sobral,  Andrews and Vacavant,  Antoine},
  year = {2014},
  month = may,
  pages = {4–21}
}

@inproceedings{brutzer2011evaluation,
  title={Evaluation of background subtraction techniques for video surveillance},
  author={Brutzer, Sebastian and H{\"o}ferlin, Benjamin and Heidemann, Gunther},
  booktitle={CVPR 2011},
  pages={1937--1944},
  year={2011},
  organization={IEEE}
}

@inproceedings{gutchess2001background,
  title={A background model initialization algorithm for video surveillance},
  author={Gutchess, Daniel and Trajkovics, M and Cohen-Solal, Eric and Lyons, Damian and Jain, Anil K},
  booktitle={Proceedings Eighth IEEE International Conference on Computer Vision. ICCV 2001},
  volume={1},
  pages={733--740},
  year={2001},
  organization={IEEE}
}

@article{Keegan2026,
  title = {Projected tensor-tensor products for efficient computation of optimal multiway data representations},
  volume = {729},
  ISSN = {0024-3795},
  url = {http://dx.doi.org/10.1016/j.laa.2025.09.018},
  DOI = {10.1016/j.laa.2025.09.018},
  journal = {Linear Algebra and its Applications},
  publisher = {Elsevier BV},
  author = {Keegan,  Katherine and Newman,  Elizabeth},
  year = {2026},
  month = jan,
  pages = {100–147}
}

@article{LizKatherine2024PROJ,
  title={Optimal matrix-mimetic tensor algebras via variable projection},
  author={Newman, Elizabeth and Keegan, Katherine},
  journal={SIAM Journal on Matrix Analysis and Applications},
  volume={46},
  number={3},
  pages={1764--1790},
  year={2025},
  publisher={SIAM}
}

@article{kong2021tensor,
  title={Tensor Q-rank: New data dependent definition of tensor rank},
  author={Kong, Hao and Lu, Canyi and Lin, Zhouchen},
  journal={Machine Learning},
  volume={110},
  number={7},
  pages={1867--1900},
  year={2021},
  publisher={Springer}
}

@article{koor2023short,
  title={A short tutorial on Wirtinger Calculus with applications in quantum information},
  author={Koor, Kelvin and Qiu, Yixian and Kwek, Leong Chuan and Rebentrost, Patrick},
  journal={arXiv preprint arXiv:2312.04858},
  year={2023}
}

\end{document}